\documentclass[a4paper]{article}
\usepackage[a4paper,top=3cm,bottom=2cm,left=3cm,right=3cm,marginparwidth=1.75cm]{geometry}
\usepackage[utf8]{inputenc} 
\usepackage[T1]{fontenc}    
\usepackage{url}           
\usepackage{booktabs}       
\usepackage{amsfonts} 
\usepackage{amsmath}
\usepackage{amssymb}
\usepackage{nicefrac}       
\usepackage[round]{natbib}
\usepackage{graphicx}
\usepackage{caption}

\newcommand{\reals}{\ensuremath{\mathbb{R}}}
\newtheorem{example2}{Example}[section]

\providecommand{\keywords}[1]
{
  \small	
  \textbf{\textit{Keywords---}} #1
}

\title{Efficient computation of matrix-vector products with full observation weighting matrices in data assimilation}

\author{Guannan Hu$^{1,*}$ and Sarah L. Dance$^{1,2}$ \\
     \small $^1$Department of Meteorology, University of Reading, United Kingdom\\
     \small $^2$Department of Mathematics and Statistics, University of Reading, United Kingdom\\
     \small $^*$Correspondence: guannan.hu@reading.ac.uk\\
}

\date{} 

\begin{document}
\maketitle

\begin{abstract}
Recent studies have demonstrated improved skill in numerical weather prediction via the use of spatially correlated observation error covariance information in data assimilation systems. In this case, the observation weighting matrices (inverse error covariance matrices) used in the assimilation may be full matrices rather than diagonal. Thus, the computation of matrix-vector products in the variational minimization problem may be very time-consuming, particularly if the parallel computation of the matrix-vector product requires a high degree of communication between processing elements. Hence, we introduce a well-known numerical approximation method, called the fast multipole method (FMM), to speed up the matrix-vector multiplications in data assimilation. We explore a particular type of FMM that uses a singular value decomposition (SVD-FMM) and adjust it to suit our new application in data assimilation. By approximating a large part of the computation of the matrix-vector product, the SVD-FMM technique greatly reduces the computational complexity compared with the standard approach. We develop a novel possible parallelization scheme of the SVD-FMM for our application, which can reduce the communication costs. We investigate the accuracy of the SVD-FMM technique in several numerical experiments: we first assess the accuracy using covariance matrices that are created using different correlation functions and lengthscales; then investigate the impact of reconditioning the covariance matrices on the accuracy; and finally examine the feasibility of the technique in the presence of missing observations. We also provide theoretical explanations for some numerical results. Our results show that the SVD-FMM technique can compute the matrix-vector product with good accuracy in a wide variety of circumstances, and hence, it has potential as an efficient technique for assimilation of a large volume of observational data within a short time interval.
\end{abstract}\hspace{10pt}


\section{Introduction}

In variational data assimilation \citep[e.g.,][]{lorenc2000met,rawlins2007met}, a nonlinear least squares problem is solved, where observations and model forecasts are blended, taking account of their uncertainties. The minimization procedure involves time-consuming computations of large matrix-vector products. In this paper we focus on the matrix-vector products involving the observations. These take the form $\mathbf{R}^{-1} \mathbf{d}$, where $\mathbf{R}^{-1} \in \reals^{m\times m}$ is the inverse of the observation error covariance matrix and $\mathbf{d} \in \reals^m$ is the observation-minus-model departure vector (see section \ref{sec:parallel standard mat-vec multiplication}).

In some practical numerical weather prediction applications, observation errors are assumed to be uncorrelated, resulting in the matrix $\mathbf{R}$ being diagonal. This reduces the number of operations required to compute matrix-vector products in the minimization, and is a pragmatic strategy when the characteristics of the observation uncertainty are not well understood \citep[e.g.,][]{liu2003potential}. However, a number of idealized and operational studies have shown that there are significant benefits to treating full observation error covariance matrices in terms of analysis information content, analysis accuracy and forecast skill, even when knowledge of the observation error correlations is only approximate \citep[e.g.,][]{HealyAndWhite2005,StewartEtAl2008,stewart2013data,weston2014accounting,simonin2019pragmatic,Bedard2020}. In particular, implementation of spatial observation error correlations modifies the length-scales of the observation increments computed by the assimilation \citep[e.g.,][]{FowlerEtAl2018,Rainwater2015,simonin2019pragmatic}, which may be especially important for multiscale systems such as convection-permitting numerical weather prediction and reanalyses \citep[e.g.,][]{Dance2019,Hu2020} or coupled ocean-atmosphere systems \citep[e.g.,][]{Hu17}. Furthermore, practical methods have been developed to estimate the observation uncertainty characteristics for a range of observation types (see the reviews by \citet{janjic2018representation} and \citet{Tandeo2020}). For example, Doppler radar radial winds, geostationary satellite data and atmospheric motion vectors (AMVs) have been shown to exhibit strong spatial error correlations \citep{waller2016DRW-UKV,waller2016SEVIRI,Cordoba17,Waller2019DRW-DWD,Honda18}. We note that in practical applications, the matrix $\mathbf{R}$ is typically treated as block diagonal with one block per observation type (for a given time). The observation errors between different types of observations are assumed to be uncorrelated.

To give an idea of the expected size of the spatially correlated observation error covariance matrices, we consider geostationary satellite observations from the SEVIRI (Spinning Enhanced Visible and Infrared Imager) instrument \citep{waller2016SEVIRI,Michel2018SEVIRI,Schmetz2002MSG}. SEVIRI measures top of the atmosphere radiances in 12 channels, with a 15-minute repeat cycle and at approximately 3 km spatial resolution  \citep[excluding the high resolution visible (HRV) channel;][]{waller2016SEVIRI,Schmetz2002MSG}. Each image consists of around $10^7$ pixels for the thermal infra-red and solar channels \citep{Schmetz2002MSG}. Hence, there are a vast number of SEVIRI observations, even for a regional domain. For instance, the number of SEVIRI observations within the domain covered by the operational AROME model from Météo-France is around $4 \times 10^5$ \citep{Seity2011MeteoFranceAROME}. In operational assimilation applications, to avoid treating spatially correlated observation errors, spatial thinning of SEVIRI observations is typically applied and this reduces the number of observations by several orders of magnitude, but also reduces their benefits in improving forecast skill \citep{waller2016SEVIRI,Michel2018SEVIRI}. The next generation of meteorological satellites will produce observations with even higher spatial resolutions \citep[][Table 1]{WMOnext}. For example, the Flexible Combined Imager (FCI) on the Meteosat Third Generation (MTG) satellite and the Advanced Baseline Imager (ABI) on the Geostationary Operational Environmental Satellite-R \citep{GOES-R} both take measurements at approximately 0.5 - 2  km spatial resolution, and the Advanced Geosynchronous Radiation Imager (AGRI) on the Fengyu-4 geostationary meteorological satellite produces observations at a resolution from  1-4 km \citep{FY-4}.

Accounting for spatially correlated error statistics in the data assimilation algorithm has potential to increase the computational cost for several reasons: (1) the computation of large matrix-vector products using parallel computing techniques, (2) the need to compute products using the inverse covariance matrix which may be different in each assimilation cycle, (3) changes to the convergence behaviour of the minimization procedure. We now review these issues in more detail. 

The first issue is that the parallel computation of a large matrix-vector product with observation data distributed across multiple processing elements (PEs) may become expensive due to excessive communications between PEs and load imbalance overheads \citep{DengBook}. Different partitions of the matrix will lead to different parallelization schemes, which will be described in section \ref{sec:parallel standard mat-vec multiplication}. In general, in order to reduce the time spent on communication operations, it is necessary to reduce the number and/or size of messages transferred. In the context of data assimilation, much of this communication can be avoided if all the observations with mutually correlated errors are assigned to one PE, as in \cite{simonin2019pragmatic}. However, this relies on the data volume of observations with mutually correlated errors being small enough for the product to be calculated on one node. It is an open question how to implement computations for larger data volumes with distributed data. While our paper is focused on parallel implementations for variational assimilation, studies such as \cite{Anderson2003} and \cite{NinoRuiz2019} address parallel implementations for ensemble methods.

The second issue is that there is a need to use the inverse observation error covariance matrix ($\mathbf{R}^{-1}$) in the computations. The most commonly used observation uncertainty diagnosis techniques provide an estimate of the observation error covariance matrix itself rather than its inverse \citep{desroziers2005diagnosis}. Since the observation distribution changes each assimilation cycle (due to quality control etc.) there is a different inverse matrix each cycle. \cite{simonin2019pragmatic} deal with this by using a Cholesky decomposition method \citep{golub1996cf} which avoids the need to compute the inverse covariance matrix directly and is applicable to any form of error covariance matrix. \cite{Guillet19} model the inverse of a spatially correlated observation error covariance matrix directly using a diffusion operator and fast unstructured meshing techniques. However, this method only deals with spatial error correlations and it is unclear if this approach is also suitable when spatial error correlations and inter-channel error correlations are combined.

The third issue is that the use of correlated observation error statistics changes the convergence behaviour of the minimization procedure. Indeed, the pre-operational experiments of \cite{weston2014accounting} showed problems with the minimization that were improved by reconditioning the observation error covariance matrix. The convergence of the minimization procedure has been further studied by \cite{tabeart2018conditioning,tabeart2020impact,tabeart2021new} who found that the minimum eigenvalue of the observation error covariance matrix is a key parameter governing the speed of convergence. Hence, the use of reconditioning methods is common in operational applications with correlated observation error statistics \citep[e.g.,][]{bormann2016enhancing,campbell2017accounting,tabeart2020impact,tabeart2020improving}.  

The aim of this paper is to carry out an investigation of a method to accelerate parallel matrix-vector products with distributed observational data, using a novel application of a well-known numerical approximation algorithm, the fast multipole method \citep[FMM;][]{rokhlin1985rapid,GREENGARD1997280}. We explore a particular type of FMM that uses a singular value decomposition (SVD) \citep{gimbutas2003generalized}. We call this method the SVD-FMM. The key idea of the SVD-FMM is to split up the computation of the matrix-vector product into separate near-field and far-field calculations. The near-field calculations are done by a standard matrix-vector multiplication of a local sub-matrix and sub-vector. The far-field calculations are carried out using the singular values and singular vectors of the sub-matrices of $\mathbf{R}^{-1}$. The sub-matrices and sub-vectors are determined by selecting specific rows and columns of $\mathbf{R}^{-1}$ corresponding to a partition of observations in the domain. The number of singular vectors employed in the calculation is chosen by the user. We will show that usually only a small number of singular vectors is needed to obtain a good accuracy (section \ref{sec:results}). The SVD-FMM technique reduces the number of floating point operations required compared with the standard approach for matrix-vector multiplication. Additionally, we develop a novel possible parallel algorithm for the SVD-FMM for our application in data assimilation, which can greatly reduce the costs of communication between PEs. The SVD-FMM allows us to assimilate a large volume of observational data in a short time interval, and has the potential to be used in practical applications.

In this initial investigation, our numerical results focus on evaluating the accuracy of the SVD-FMM. Thus our experiments are carried out for an idealized problem, using serial, rather than parallel computing. Nevertheless, we do compare the efficiency of the proposed parallelization scheme with different parallel formulations of standard matrix-vector multiplications in terms of communication costs. In our current implementation, the method needs to be applied after we obtain a representation of the matrix $\mathbf{R}^{-1}$. In order to have a clear focus solely on computing matrix-vector products we do not address the computation of the inverse observation error covariance matrix. Instead, we assume that the inverse observation error covariance is already known. 

In our experiments we show that the SVD-FMM can work well with the inverses of a variety of covariance matrices. In particular, we apply the SVD-FMM to the inverses of covariance matrices created using different correlation functions and lengthscales. We also use the reconditioning methods of \cite{tabeart2020improving} together with the SVD-FMM to gain insight into how the accuracy of the SVD-FMM should change with different levels of reconditioning. In practice, the observation distribution varies each assimilation cycle due to factors such as quality control and the removal of cloudy satellite radiance observations. Therefore, we further carry out some experiments to demonstrate that the SVD-FMM is feasible even if there are some missing observations.

The rest of this paper is organised as follows: We provide the parallel formulations of standard matrix-vector multiplication in section \ref{sec:parallel standard mat-vec multiplication}. In section \ref{sec:application of FMM to DA} we present our novel algorithm for the SVD-FMM and compare the complexity and the communication costs between the SVD-FMM and standard, parallel matrix-vector multiplication. In section \ref{sec:experimental design} we explain our experimental design for our idealized experiments. In particular, we describe the generation of observations and observation error covariance matrices as well as the reconditioning methods. In section \ref{sec:results} we show the results of our numerical experiments with varying correlation functions, length scales, reconditioning methods and condition numbers. We also show how the results change with missing observations. Finally we give a summary in section \ref{sec:summary} and conclude that our proposed algorithm has potential for use in operational data assimilation for fast computation of large matrix-vector products.

\section{Parallelization of standard matrix-vector multiplication}\label{sec:parallel standard mat-vec multiplication}

In this section we describe three distinct standard parallel formulations for computing large matrix-vector products \citetext{\citealp[see][Section 8.1]{grama2003introduction}; \citealp[][Section 6.2.1]{DengBook}}. We will also discuss how to exploit the symmetric structure of $\mathbf{R}^{-1}$ for parallelization \citep{Geus2001Symmetry}. These matrix-vector products arise in the solution of the variational minimization problem in data assimilation in the form
\begin{equation} \label{eq:mat-vec product}
	\mathbf{q}=\mathbf{A}\mathbf{d},
\end{equation}
where $\mathbf{A} \in \mathbb{R}^{m \times m}$ denotes the inverse of the observation error covariance matrix, $\mathbf{d} \in \mathbb{R}^m$ denotes the observation-minus-model departure vector, and $\mathbf{q} \in \mathbb{R}^m$ denotes the result of the multiplication. For the purposes of this description, we assume that the matrix $\mathbf{A}$ and vector $\mathbf{d}$ are already known. 

To see how matrix-vector products in the form of Eq. \eqref{eq:mat-vec product} arise in data assimilation, we consider the observation penalty term of the variational data assimilation cost function. This is given by multiplying the inverse observation error covariance matrix by the observation-minus-model differences 
\begin{equation}\label{eq:observation penalty}
    J_o=\frac{1}{2}\{\mathbf{y}-H(\mathbf{x})\}^\mathrm{T}\mathbf{R}^{-1}\{\mathbf{y}-H(\mathbf{x})\},
\end{equation}
where $\mathbf{y} \in \reals^m$ denotes the observation vector, $\mathbf{x} \in \reals^n$ denotes model state and $H$ denotes the observation operator that maps model state to the observation \citep[$H:\reals^n \rightarrow \reals^m$; e.g.,][]{lorenc2000met,rawlins2007met,simonin2019pragmatic,Nichols2010}. To solve the variational minimization problem, the gradient of Eq. \eqref{eq:observation penalty} is needed \citep[e.g.,][]{simonin2019pragmatic}
\begin{equation}
    \frac{\partial J_o}{\partial \mathbf{x}}=\mathbf{H}^\mathrm{T}\mathbf{R}^{-1}\{\mathbf{y}-H(\mathbf{x})\}.
\end{equation}
Thus matrix-vector products of the form of Eq. \eqref{eq:mat-vec product} arise in both the computation of the cost function and its gradient, with $\mathbf{A} = \mathbf{R}^{-1}$ and $\mathbf{d}=\mathbf{y}-H(\mathbf{x})$.

\begin{figure}[!ht]
    \centering
    \includegraphics[width=.75\textwidth]{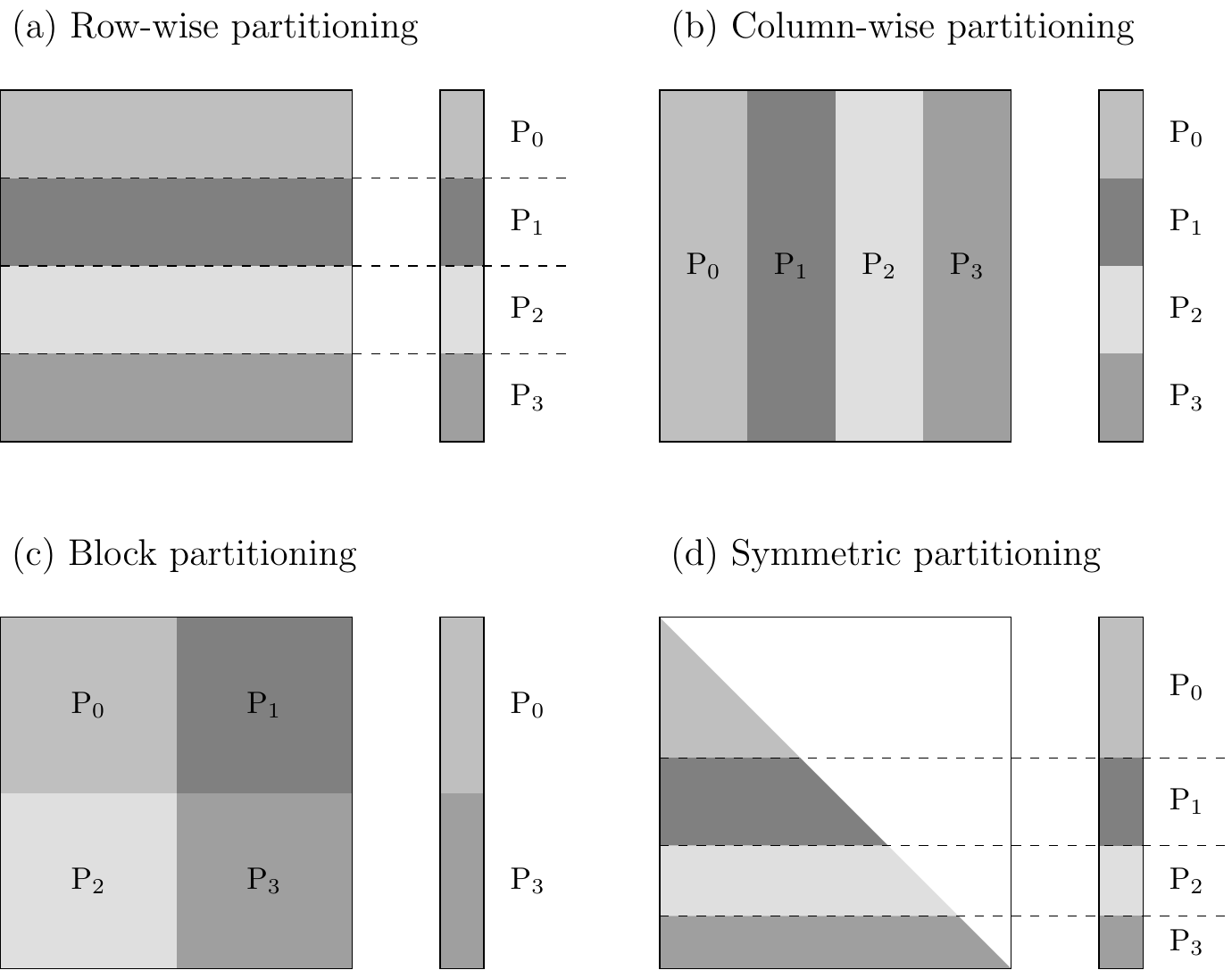}
	\caption{Different ways of partitioning matrices (represented by squares) and vectors (represented by bars) for parallel computation of the matrix-vector products. Different colours demonstrate the portions of matrices and vectors that are distributed over four PEs ($P_0$, $P_1$, $P_2$ and $P_3$).} \label{fig:matrix partitioning}
\end{figure}

The parallelization schemes for computing Eq. \eqref{eq:mat-vec product} start with the distribution of observations over PEs. We consider a simple domain decomposition, in which the observations are distributed over a number of PEs according to their geophysical locations. This will result in a split of the components of matrix $\mathbf{A}$ and vector $\mathbf{d}$ across PEs. We will introduce four different partitions of the matrix (see Fig. \ref{fig:matrix partitioning}). Each of them will lead to a unique parallelization scheme. The communication costs of each scheme depend on various parameters, including i) the time to prepare a message for transmission, ii) the time it takes for a message to travel (latency), iii) how many words can traverse per second (bandwidth), iv) how many PEs to communicate with and v) the message size \citep{grama2003introduction}. Since the first three parameters are determined by the configuration of the parallel machine, we discuss the communication costs for different parallelization schemes using the last two parameters.

\subsection{Row-wise partitioning}\label{sec:row-wise}

We first consider a row-wise partitioning, in which case each PE stores one row or several rows of $\mathbf{A}$ and one element or a portion of $\mathbf{d}$. Fig. \ref{fig:matrix partitioning}(a) illustrates the partitioning using four PEs. In order for each PE to perform its computation, we need to distribute the full vector among all the PEs. This requires an all-to-all broadcast. Based on table 4.1 in \cite{grama2003introduction}, the communication cost of this communication operation is $t_s\log B +t_wm$, where $B$ is the number of PEs, $t_s$ is the startup time and $t_w$ is the per-word transfer time \citep{grama2003introduction}. The two time parameters depend on computer architecture and performance. After the communication, each PE multiplies its $m/B$ rows with the vector, which requires on the order of $m^2/B$ floating point operations.

\subsection{Column-wise partitioning}\label{sec:column-wise}

Instead of storing row(s), each PE can store the column(s) of $\mathbf{A}$. Fig. \ref{fig:matrix partitioning}(b) shows an example of using four PEs. With a column-wise partitioning, each PE can calculate a partial result locally. Each PE multiplies $m/B$ columns of $\mathbf{A}$ with $m/B$ elements of $\mathbf{d}$, which requires $\mathcal{O}(m^2/B)$ floating point operations. After the local computation, we need to perform an all-to-one reduction to sum up the partial results given by each PE, which takes time $(t_s+t_wm)\log B$ \citep[][Table 4.1]{grama2003introduction}. After the all-to-one reduction only one PE contains the full vector $\mathbf{q}$ and thus, we may need to redistribute $\mathbf{q}$. This requires a scatter operation which takes time $t_s\log B+t_wm$ \citep[][Table 4.1]{grama2003introduction}. Although the column-wise partitioning circumvents the use of an expensive all-to-all broadcast operation, it increases the message size for data-transfer operations. The message size for the all-to-one reduction is $m$, whereas the message size for the all-to-all broadcast used for the row-wise partitioning is $m/B$. It should be noticed that in data assimilation applications, the matrix $\mathbf{A}$ is known to be symmetric, which means that the parallel algorithm using column-wise partitioning can also be used with row-wise partitioning.

\subsection{Block 2-D partitioning}\label{block 2D}

The row-wise partitioning and column-wise partitioning are 1-D partitioning. We now consider a 2-D partitioning, which distributes the blocks of $\mathbf{A}$ among PEs. An example of the block 2-D partitioning is shown in Fig. \ref{fig:matrix partitioning}(c). Note that the matrix $\mathbf{A}$ is equally separated by $4$ PEs, while the vector $\mathbf{d}$ is only distributed among $2$ PEs that own the diagonal blocks of $\mathbf{A}$, i.e., $\text{P}_0$ and $\text{P}_3$. The first step is for each PE that stores a portion of $\mathbf{d}$ to broadcast that portion to the other PEs in the same column, which requires a column-wise one-to-all broadcast. The communication time of this operation is $(t_s+t_wm/\sqrt{B})\log \sqrt{B}$. The next step is to perform a row-wise all-to-one reduction, which requires another $(t_s+t_wm/\sqrt{B})\log \sqrt{B}$ times. \cite{grama2003introduction} showed that computing the matrix-vector product using the block 2-D partitioning of the matrix can be faster than using 1-D partitioning for the same number of PEs. However, a potential problem is that when using the block 2-D partitioning, the vector elements are not distributed among all the PEs and this can result in load imbalance.

\subsection{Partitioning for symmetric matrices}\label{sec:symmetric}

A possible partitioning of $\mathbf{A}$ and $\mathbf{d}$ taking into account the symmetric structure of $\mathbf{A}$ is shown in Fig. \ref{fig:matrix partitioning}(d). The first step is to exchange the portions of $\mathbf{d}$ among PEs. This requires an all-to-all broadcast with an average message size of $m/B$. It should be noticed that not all PEs need to send their portion of $\mathbf{d}$ to all other PEs. In the case of four PEs, $P_0$ sends data to $P_1$, $P_2$ and $P_3$, $P_1$ to $P_2$ and $P_3$, and $P_2$ to $P_3$. The second step is to multiply the vector elements with the local part of the matrix. The last step is to exchange the results of local computation using an all-to-one reduction. The average message size for this operation is $m/B$. In the case of four PEs, $P_0$ collects the results from all other PEs, $P_1$ collects the results from $P_2$ and $P_3$, and $P_2$ collects the result from $P_3$. The last PE does not need to collect the results from other PEs. Slightly different parallelization schemes may be used for this kind of partitioning, depending on the computer architecture \citep{Geus2001Symmetry}. This partitioning saves the storage space of each PE (only half of the matrix is stored) and keeps the load balanced. Furthermore, this partitioning has a very large advantage for sparse matrices, in which case each PE only needs to communicate with neighbouring PEs \citep{Geus2001Symmetry}.

In the next section, we describe a different approach for calculating large matrix-vector products, that is applicable to any matrix (most useful for full or dense matrices). This approach reduces the computational cost for serial calculations of Eq. \eqref{eq:mat-vec product} and the message size for communication operations.\\

\begin{figure}[!ht]
	\centering
	\includegraphics[width=.6\textwidth]{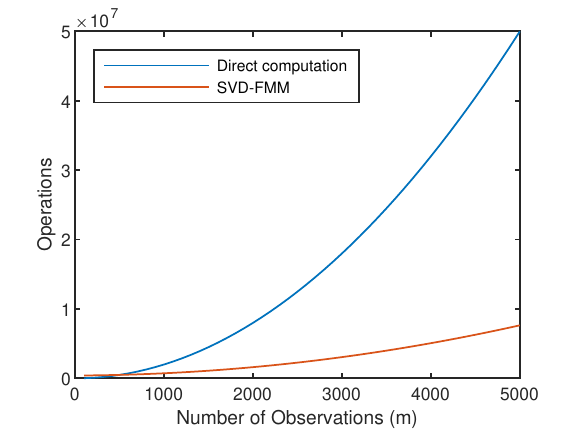}
	\caption{The number of floating point operations for the direct computation of the matrix-vector product and SVD-FMM with our configuration of boxes and $p=10$.}
	\label{fig:complexity}
\end{figure}

\section{Application of the fast multipole method (FMM) to data assimilation} \label{sec:application of FMM to DA}

The fast multipole method (FMM) was originally presented for the rapid evaluation of all pairwise interactions between a large number of charged particles \citep{rokhlin1985rapid,GREENGARD1997280}. The mathematical form of this fast summation is equivalent to Eq. \eqref{eq:mat-vec product}. While the direct calculation of Eq. \eqref{eq:mat-vec product} requires $\mathcal{O}(m^2)$ work, the FMM needs only $\mathcal{O}(m)$ operations. In addition, the FMM relies on a hierarchical division of the computational domain, which is well suited for parallel computing. In the classical FMM, the matrix $\mathbf{A}$ is given by some functions describing pairwise interactions between particles, such as those describing gravitational and electrostatic potentials. In our application, however, the matrix $\mathbf{A}$ is obtained by inverting the matrix $\mathbf{R}$. Therefore, we need a generalized FMM that does not require an underlying analytic function \citep{gimbutas2003generalized}.

\subsection{Separation of the matrix-vector product into near- and far-field terms} \label{sec:FMM approach}

Since the SVD-FMM is not commonly used in meteorology, we first explain the key idea behind the technique and give a simple example to aid the reader's understanding, before going into the mathematical details in later subsections. The basic idea of the SVD-FMM is to calculate Eq. \eqref{eq:mat-vec product} in two parts:
\begin{equation} \label{eq:near and far-field mat-vec product}
	\mathbf{q}=\mathbf{A}(\cdot,\mathbf{I}_1)\mathbf{d}(\mathbf{I}_1) + \mathbf{A}(\cdot,\mathbf{I}_2)\mathbf{d}(\mathbf{I}_2),
\end{equation}
where the column indices of $\mathbf{A}$ and the indices of $\mathbf{d}$ are divided into two mutually independent sets denoted by $\mathbf{I}_1$ and $\mathbf{I}_2$. The indices are selected according to a partition of observations (described in more detail in section \ref{sec:box-tree}). The first part of Eq. \eqref{eq:near and far-field mat-vec product}, referred to as the near-field calculation, is computed using a standard matrix-vector multiplication and the second part, referred to as the far-field calculation, is estimated using an approximation. In the SVD-FMM, the far-field calculation is computed using the singular value decomposition (SVD) of $\mathbf{A}(\cdot,\mathbf{I}_2)$. We provide a simple example to aid the reader's understanding.

\begin{example2}
    Suppose we have a matrix $\mathbf{A} \in \reals^{2 \times 3}$ and a vector $\mathbf{d} \in \reals^3$ given by
    \begin{align*}
        \mathbf{A} = \begin{pmatrix}
            a_{11} & a_{12} & a_{13}  \\
            a_{21} & a_{22} & a_{23}     
        \end{pmatrix} \quad \text{and} \quad
        \mathbf{d} = \begin{pmatrix}
            d_{1} \\
            d_{2} \\
            d_{3} \\
        \end{pmatrix}.
    \end{align*}
    A standard matrix-vector multiplication gives 
    \begin{align*}
        \mathbf{q} = \mathbf{A}\mathbf{d}=  \begin{pmatrix}
        a_{11}d_{1} + a_{12}d_{2} + a_{13}d_{3} \\
        a_{21}d_{1} + a_{22}d_{2} + a_{23}d_{3} 
        \end{pmatrix}.
    \end{align*}
    We can also compute $\mathbf{q}$ by partitioning the problem as in Eq. \eqref{eq:near and far-field mat-vec product} by letting $\mathbf{I}_1=\{2\}$ (the near field) and $\mathbf{I}_2=\{1,3\}$ (the far field), such that
    \begin{align} \label{eq:example near and far field}
        \mathbf{q}= \begin{pmatrix}
            a_{12} \\
            a_{22}  
        \end{pmatrix} \cdot d_{2} + 
        \begin{pmatrix}
            a_{11} & a_{13} \\
            a_{21} & a_{23}
        \end{pmatrix}
        \cdot \begin{pmatrix}
            d_{1}  \\
            d_{3}   
        \end{pmatrix}. 
    \end{align}
    The singular value decomposition (SVD) of the sub-matrix in the second term is given by
    \begin{align*}
        \begin{pmatrix}
            a_{11} & a_{13} \\
            a_{21} & a_{23}
        \end{pmatrix} =  \sum_{i=1}^2 
        \mathbf{u}_i \sigma_i (\mathbf{v}_i)^T,
    \end{align*}
    where $\mathbf{u}_1, \mathbf{u}_2 \in \reals^2$ are the orthonormal left singular vectors, $\sigma_1, \sigma_2$ are the scalar singular values, and  $\mathbf{v}_1, \mathbf{v}_2 \in \reals^2$ are the orthonormal right singular vectors. If $\sigma_1 \gg \sigma_2$ then we can truncate the SVD to give an approximation for the far field term of Eq. \eqref{eq:example near and far field} as
    \begin{align*}
        \mathbf{q} \approx \begin{pmatrix}
            a_{12} \\
            a_{22}  
        \end{pmatrix} \cdot d_{2} + \mathbf{u}_1 \sigma_1 (\mathbf{v}_1)^T \cdot \begin{pmatrix}
            d_{1}  \\
            d_{3}   
        \end{pmatrix}.
    \end{align*}
\end{example2}

More generally, for a larger problem, we can use only the first few singular vectors and singular values to estimate the far-field calculation. If the matrix $\mathbf{A}$ remains fixed, we can use the same singular vectors and singular values to compute $\mathbf{q}$ for any vector $\mathbf{d}$. This will reduce the cost for calculating Eq. \eqref{eq:mat-vec product} (unless the matrix $\mathbf{A}$ is sparse). Furthermore, storing the singular vectors and singular values requires less memory than storing the full far-field sub-matrix. Note that if the matrix $\mathbf{A}$ changes too often, then we would need to perform the SVD again and again, which would be computationally expensive.

The example we have discussed in this section only shows the basic idea of the SVD-FMM. The actual SVD-FMM algorithm is multi-level and more complicated. Before we can address this, we explain the partitioning of the observations that is needed for the multi-level algorithm.

\subsection{Partition of observations into a hierarchical structure of nested boxes}\label{sec:box-tree}

\begin{figure}
	\centering
    \includegraphics{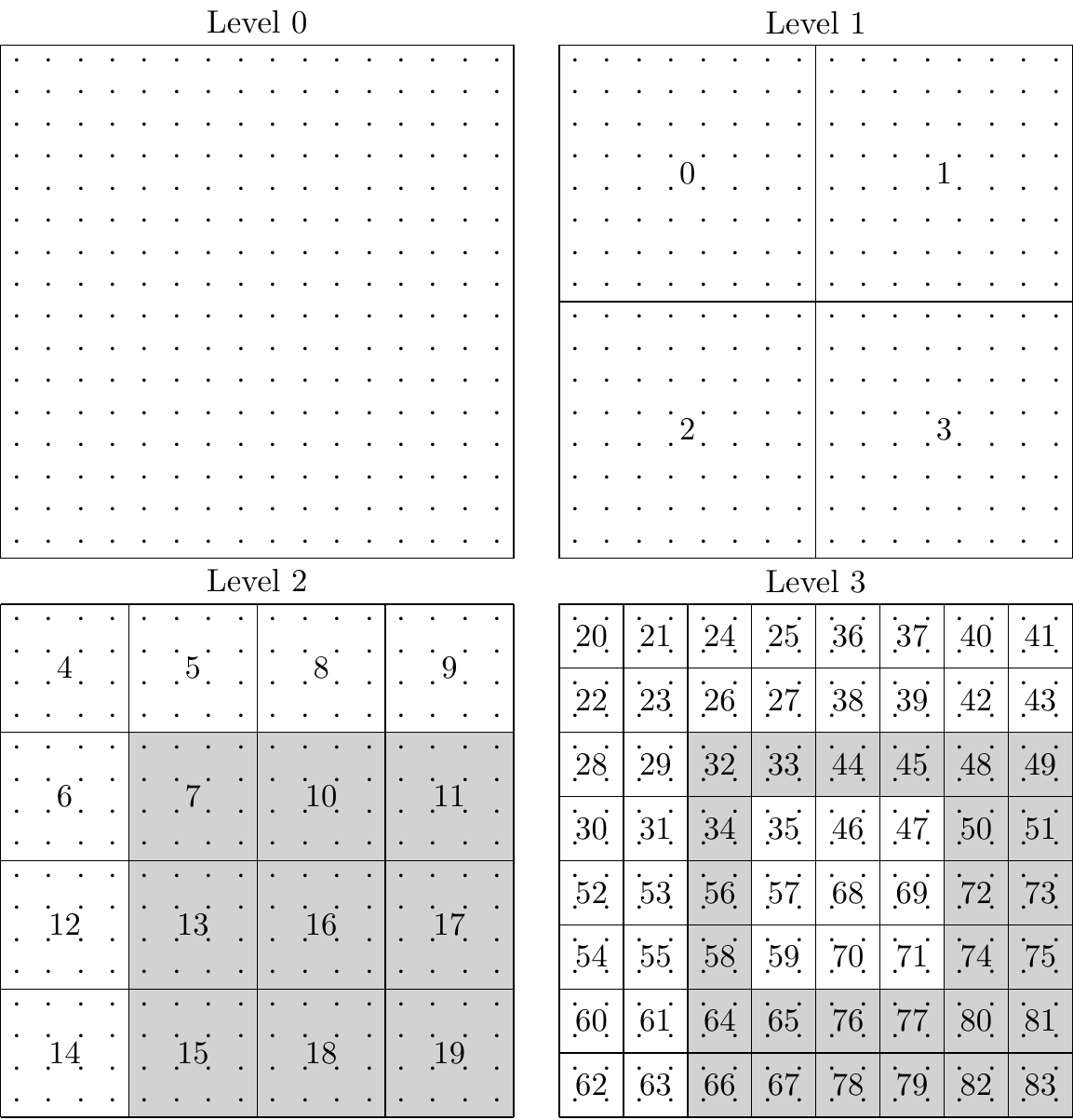}   
	\caption{Illustration of the hierarchy of boxes of level 0-3 in our example. Level 0 refers to the box that covers the entire observation domain, while level $l+1$ refers to the boxes that are obtained by equally subdividing each box in level $l$ into four. The dots represent the observation locations and the numbers are the box indices. The shaded boxes in level 2 are the boxes that are in the near-field of box 16, while the shaded boxes in level 3 are the boxes that are in the interaction list of box 68. The definitions of near-field and interaction list can be found in the main text.}\label{fig:division}
\end{figure}

In this section, we describe the partition of observations and explain the notation that we use. Note that we use the nomenclature found in the FMM literature \citep[e.g.,][]{gimbutas2003generalized}. Suppose that we have $m$ observations. These could be located on a latitude-longitude grid, or irregularly distributed. We first find a minimal square (or rectangle) that covers the locations of all of the observations, and then hierarchically subdivide the observation domain into smaller and smaller boxes to generate a \textit{box-tree} - a hierarchical structure of nested boxes. An example is given in Fig. \ref{fig:division}. In the tree structure, \textit{level 0} refers to the biggest box that covers the entire observation domain, and \textit{level} $l+1$ refers to the boxes that are obtained from \textit{level} $l$ by subdividing each box into four smaller boxes of equal size. In our particular example, we choose level 3 as the highest level, which is the minimal number of levels required for the present SVD-FMM approach. The boxes at the highest level are called \textit{leaf boxes}.

In practical applications, the number of levels is determined such that the averaged number of observations in the leaf boxes is smaller than a prescribed value. For unevenly distributed observations, an adaptive tree structure could be created \citep{gimbutas2003generalized}, in which smaller boxes are generated only where data is dense. 
 
We number the boxes in all levels except level 0 by a \textit{Z-order curve} as shown in Fig. \ref{fig:division} \citep{gargantini1982effective}. This facilitates easy formulae for box indexing as will be described later in this subsection. We call boxes \textit{neighbours} if they are on the same level and connect to each other. In each level the \textit{near-field} of a box $b$ is made of itself and all its neighbours.  The \textit{far-field} of the box $b$ is made of everything else. We use $\mathcal{N}_b$ and $\mathcal{F}_b$ to denote the lists of boxes that are in the near-field and far-field of box $b$, respectively. For example, the near-field of box 16 in Fig. \ref{fig:division} consists of the 9 shaded boxes, namely, $\mathcal{N}_{16} = \{7, 10, 11, 13, 15, 16, 17, 18, 19\}$, and the far-field of box 16 is the other 7 boxes, i.e., $\mathcal{F}_{16} = \{4, 5, 6, 8, 9, 12, 14\}$. The smallest number of boxes in a box's near-field is 4, which occurs when the box is on a corner of the observation domain. 

In a box-tree the \textit{children} of box $b$ in level $l$, denoted by $\mathcal{C}_b$, refer to the four boxes in level $l+1$ that are subdivided from itself, such as $\mathcal{C}_4 = \{20, 21, 22, 23\}$. In the Z-curve ordering, the indices of the children of box $b$ are given by $4b+4$, $4b+5$, $4b+6$, and $4b+7$. The \textit{parent} of box $b$, denoted by $\mathcal{P}_b$, refers to the box in a coarser level that contains box $b$. For instance, $\mathcal{P}_{20} = \mathcal{P}_{21} = \mathcal{P}_{22} =\mathcal{P}_{23} = \{4\}$. 
 
The \textit{interaction list} of box $b$, denoted by $\mathcal{L}_b$, is the set of boxes which are children of the neighbours of box $b$'s parents and which are in the far-field of box $b$. For example, $\mathcal{L}_{68}=\{32,33,34,44,45,48,49,50,51,56,\allowbreak58,64,65,66,67,72,73,74,75,76,77,78,79,80,81,82,83\}$ and $\mathcal{L}_{16}=\{4,5,6,8,9,12,14\}$. Note that the interaction list of a box in level 2 is exactly the same as its far-field. 

\subsection{The SVD-FMM algorithm} \label{sec:SVD-FMM}

In this section we describe the multi-level SVD-FMM algorithm \citep{gimbutas2003generalized}. For concreteness we choose a particular configuration of boxes (as shown in Fig. \ref{fig:division}) to describe the algorithm, but it can be generalized to other configurations. For each box $b$ we can write Eq. \eqref{eq:near and far-field mat-vec product} as
\begin{equation} \label{eq:near and far-field mat-vec product box}
	\mathbf{q}(\mathbf{I}_b)=\mathbf{A}(\mathbf{I}_b,\mathbf{I}_{\mathcal{N}_b})\mathbf{d}(\mathbf{I}_{\mathcal{N}_b}) + \mathbf{A}(\mathbf{I}_b,\mathbf{I}_{\mathcal{F}_b})\mathbf{d}(\mathbf{I}_{\mathcal{F}_b}),
\end{equation}
where $\mathbf{I}_b$, $\mathbf{I}_{\mathcal{N}_b}$ and $\mathbf{I}_{\mathcal{F}_b}$ denote the sets of observation indices in box $b$, $\mathcal{N}_b$ and $\mathcal{F}_b$, respectively, and $\mathbf{A}(\mathbf{I}_b,\mathbf{I}_{\mathcal{N}_b}) =\{a_{i,j}|i \in \mathbf{I}_{b}, \; j \in \mathbf{I}_{\mathcal{N}_b}\}$ and $\mathbf{A}(\mathbf{I}_{b},\mathbf{I}_{\mathcal{F}_b})=\{a_{i,j}|i \in \mathbf{I}_b, \; j \in \mathbf{I}_{\mathcal{F}_b}\}$ denote sub-matrices of $\mathbf{A}$ that are comprised of specific rows and columns of $\mathbf{A}$. We call $\mathbf{q}(\mathbf{I}_b)=\{q_i|i \in \mathbf{I}_b\}$ the \textit{target} and $\mathbf{d}(\mathbf{I}_{\mathcal{F}_b})=\{d_j|j \in \mathbf{I}_{\mathcal{F}_b}\}$ the \textit{source}.

The first matrix-vector product in Eq. \eqref{eq:near and far-field mat-vec product box} is calculated  from the near-field components. It is computed directly without using any approximation. The second matrix-vector product in Eq. \eqref{eq:near and far-field mat-vec product box}, containing the far-field components, is estimated by matrix decomposition, using a multi-level approach exploiting the hierarchical box-tree structure established in section \ref{sec:box-tree}. Before giving the mathematical details, we first give an overview of the algorithm. 

\begin{figure}
    \centering
    \includegraphics{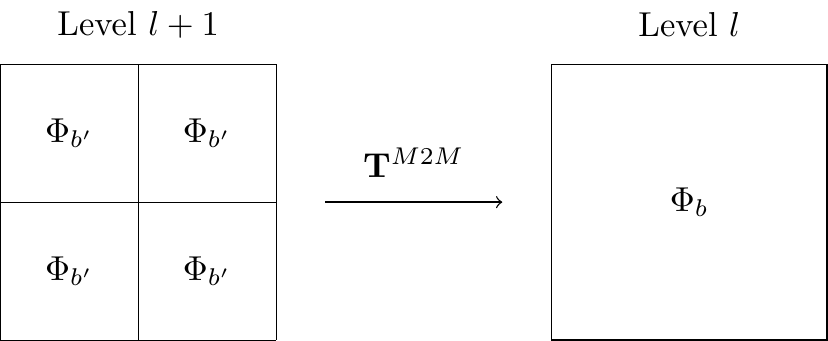}    
	\caption{Illustration of the multipole-to-multipole translation operator: $\mathbf{T}^{M2M}$ transforms the multipole expansions ($\Phi_{b'}$) of box $b'$ to the multipole expansion ($\Phi_{b}$) of box $b$, where box $b$ is the parent of box $b'$.}\label{fig:M2M}
\end{figure}

\begin{figure}
    \centering
    \includegraphics{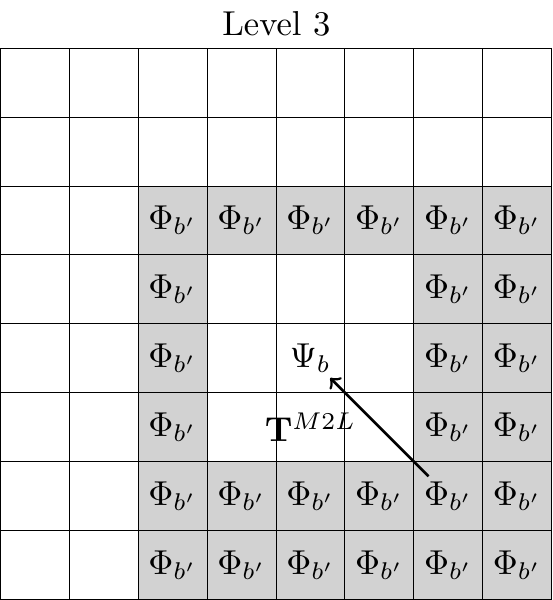}
	\caption{Illustration of the multipole-to-local translation operator: $\mathbf{T}^{M2L}$ transforms the multipole expansions ($\Phi_{b'}$) of box $b'$ to the local expansion ($\Psi_{b}$) of box $b$, where box $b'$ is in the interaction list of box $b$. The arrow illustrates the transformation from one particular $b'$ to $b$. For each $b'$ we have a unique $\mathbf{T}^{M2L}$.} \label{fig:M2L}
\end{figure}

\begin{figure}
    \centering
	\includegraphics{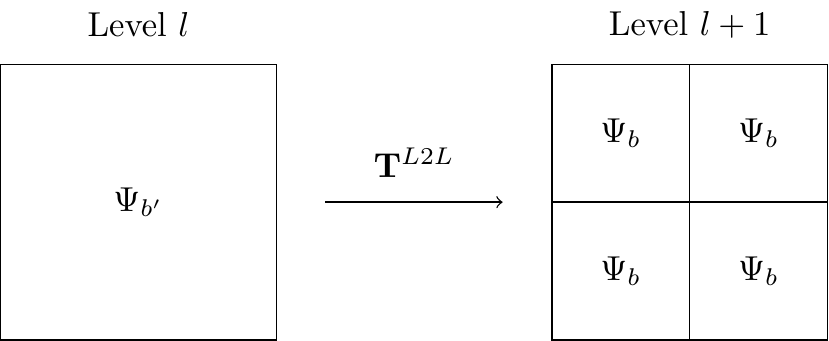}
	\caption{Illustration of the local-to-local translation operator: $\mathbf{T}^{L2L}$ transforms the local expansion ($\Psi_{b'}$) of box $b'$ to the local expansions ($\Psi_{b}$) of box $b$, where box $b$ is the child of box $b'$.}\label{fig:L2L}
\end{figure}

The key idea of the multi-level approach is to perform the SVD of the sub-matrices of $\mathbf{A}$ given by $\mathbf{A}(\mathbf{I}_b,\mathbf{I}_{\mathcal{F}_b})$ and use the singular vectors and singular values just obtained to create a \textit{multipole expansion}, a \textit{local expansion} and three \textit{translation operators}, and then use them to estimate the result. The multipole expansion can be interpreted as a short representation of the sub-vector of $\mathbf{d}$ given by $\mathbf{d}(\mathbf{I}_b)$, which is obtained by projecting the components of $\mathbf{d}(\mathbf{I}_b)$ onto the basis given by $p$ singular vectors and hence, is a short vector with $p$ components. The local expansion can be considered as the short representation of the sub-vector of $\mathbf{d}$ given by $\mathbf{d}(\mathbf{I}_{\mathcal{F}_b})$, which is also a vector with $p$ components. It is computed from the multipole expansions of a group of boxes that are in $b$'s interaction list using the translation operators. The translation operators allow exploitation of the hierarchical structure established in the box-tree, by transforming the projection from one basis of singular vectors into another basis of singular vectors. There are three kinds of translation operators: the \textit{multipole-to-multipole} (M2M) translation operator transforms the multipole expansion of a box to the multipole expansion of its parent (see Fig. \ref{fig:M2M}); the \textit{multipole-to-local} (M2L) translation operator translates the multipole expansion of a box to the local expansion of another box in the same level (see Fig. \ref{fig:M2L}); and the \textit{local-to-local} (L2L) translation operator converts the local expansion of a non-leaf box to the local expansion of its children (see Fig. \ref{fig:L2L}). 

Due to the use of these expansions and translation operators, the SVD-FMM computes matrix-vector products more efficiently than the standard approach, which is reflected in both serial efficiency (algorithmic complexity) and parallel efficiency (see sections \ref{sec:complexity} and \ref{sec:parallelism} for more details). 

\begin{figure}
    \centering
	\includegraphics{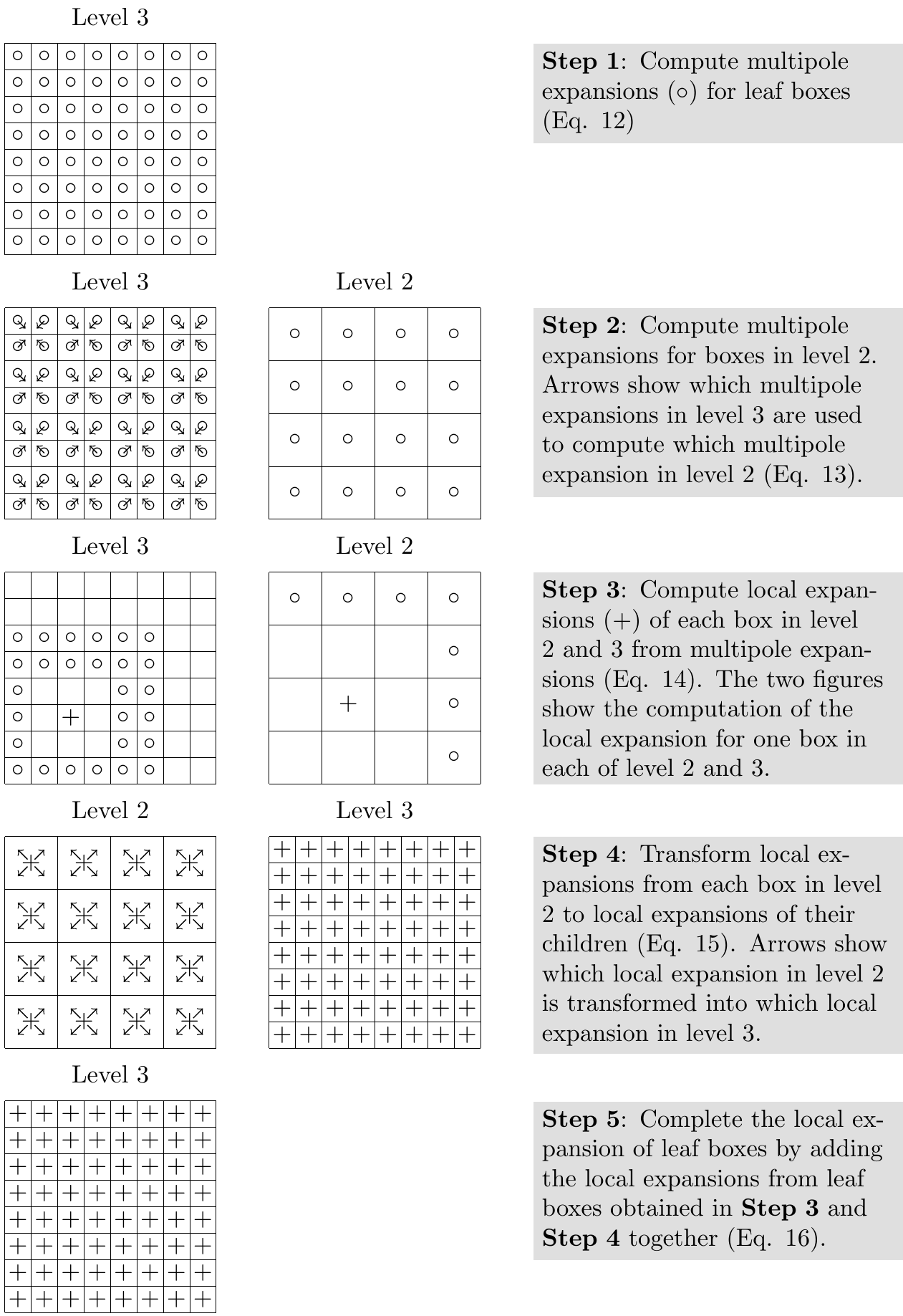}
    \caption{Illustration of the SVD-FMM algorithm.}
    \label{fig:illustration of algorithm}
\end{figure}

We now describe the mathematical steps of the SVD-FMM algorithm in more detail. A schematic illustrating the algorithm is given in Fig.~\ref{fig:illustration of algorithm}.\\

\noindent \textbf{Initialization:} \textit{Compute the SVD of sub-matrices of $\mathbf{A}$ and the translation operators}. We compute truncated SVDs of sub-matrices of $\mathbf{A}$ for each box $b$ in level 2 and 3
\begin{equation}\label{eq:decomposition}
    \begin{aligned}
        \mathbf{A}(\mathbf{I}_{\mathcal{F}_b},\mathbf{I}_b) &\approx 
        \sum_{k=1}^p\mathbf{u}_{k}^{src,b} s_{k}^{src,b} (\mathbf{v}_{k}^{src,b})^\mathrm{T}, \\ \mathbf{A}(\mathbf{I}_b, \mathbf{I}_{\mathcal{F}_b}) &\approx
        \sum_{k=1}^p\mathbf{u}_{k}^{tgt,b} s_{k}^{tgt,b} (\mathbf{v}_{k}^{tgt,b})^\mathrm{T}, 
    \end{aligned}
\end{equation}
where $\mathbf{u}^{src,b}_k  \in \mathbb{R}^{m_{\mathcal{F}_b}}$ and $\mathbf{u}^{tgt,b}_k  \in \mathbb{R}^{m_b}$ are the $k$-th left singular vectors, $s^{src,b}_k$ and $s^{tgt,b}_k$ are the $k$-th singular values, and $\mathbf{v}^{src,b}_k \in \mathbb{R}^{m_b}$ and $\mathbf{v}^{tgt,b}_k \in \mathbb{R}^{m_{\mathcal{F}_b}}$ are the $k$-th right singular vectors. The symbols $m_{\mathcal{F}_b}$ and $m_b$ denote the length of $\mathbf{I}_{\mathcal{F}_b}$ and $\mathbf{I}_b$, respectively. The symbol $p$ denotes the number of singular vectors and the number of singular values. 

(Note that matrix $\mathbf{A}$ is symmetric for the data assimilation problem discussed here. Therefore, in practice we actually only need to compute one SVD per box $b$, and for each box we have $\mathbf{u}^{tgt,b}_k=\mathbf{v}^{src,b}_k$, $\mathbf{s}^{tgt,b}_k=\mathbf{s}^{src,b}_k$ and $\mathbf{v}^{tgt,b}_k=\mathbf{u}^{src,b}_k$. However, we retain the different notations for the purpose of clarifying the algorithm.)

The singular vectors and singular values are then used to generate multipole expansions and local expansions and a set of translation operators. The M2M translation operator (see Fig.~\ref{fig:M2M}) converts the multipole expansion of a child box into the multipole expansion of its parent and is defined as 
\begin{equation}\label{eq:M2M}
    \mathbf{T}^{M2M}(k,k') = \sum_{i \in \mathbf{I}_{b'}}v^{src,b}_{k,i}v^{src,b'}_{k',i}, \quad k,k' = 1,\cdots,p
\end{equation}
for each $b' \in \mathcal{C}_b$, where $v^{src,b}_{k,i}$ and $v^{src,b'}_{k,i}$ denote the elements of $\mathbf{v}_{k}^{src,b}$ and $\mathbf{v}_{k}^{src,b'}$ that correspond to the $i$-th observation in box $b'$. The M2L translation operator (see Fig.~\ref{fig:M2L}) converts the multipole expansion of box $b'$ into the local expansion of box $b$, for each $b'$ in the interaction list of $b$. It is defined as 
\begin{equation}\label{eq:M2L}
    \mathbf{T}^{M2L}(k,k')= \sum_{i \in \mathbf{I}_{b'}} v^{tgt,b}_{k,i}v^{src,b'}_{k',i},
    \quad k,k' = 1,\cdots,p
\end{equation}
for each $b' \in \mathcal{L}_b$, where $v^{tgt,b}_{k,i}$ and $v^{src,b'}_{k,i}$ denote the elements of $\mathbf{v}_{k}^{tgt,b}$ and $\mathbf{v}_{k}^{src,b'}$ that correspond to the $i$-th observation in box $b'$. The L2L translation operator (see Fig.~\ref{fig:L2L}) transfers the local expansion of a parent box to its children and is defined as
\begin{equation}\label{eq:L2L}
    \mathbf{T}^{L2L}(k,k')=  \sum_{i \in \mathbf{I}_{\mathcal{F}_{b'}}}v^{tgt,b}_{k,i}v^{tgt,b'}_{k',i},
    \quad k,k' = 1,\cdots,p
\end{equation}
for each $b' = \mathcal{P}_b$, where $v^{tgt,b}_{k,i}$ and $v^{tgt,b'}_{k',i}$ denote the elements of $\mathbf{v}_{k}^{tgt,b}$ and $\mathbf{v}_{k}^{tgt,b'}$ that correspond to the $i$-th observation in the far-field of box $b'$.\\

\noindent \textbf{Step 1:} \textit{Compute the multipole expansion for the leaf boxes}. The multipole expansion for each leaf box is computed by
\begin{equation}\label{eq:compute multipole childless level}
    \Phi^{b}_k = \sum_{j \in \mathbf{I}_{b}}v_{k,j}^{src,b}d_j,
\end{equation}
where $b=20,\cdots,83$ and $k=1,\cdots,p$.\\

\noindent \textbf{Step 2:} \textit{Compute the multipole expansion for the non-leaf boxes}. The multipole expansion for each box in level 2 is computed by translating the multipole expansions from its children in level 3 using the $\mathbf{T}^{M2M}$ translation operator,
\begin{equation}\label{eq:compute multipole parent level}
    \Phi^{b}_k = \sum_{b' \in \mathcal{C}_b}\sum_{k'=1}^{p}\mathbf{T}^{M2M}(k,k')\Phi^{b'}_{k'},
\end{equation}
where $b=4,\cdots,19$, $k=1,\cdots,p$ and $b'$ denotes the child of $b$.\\

\noindent \textbf{Step 3:} \textit{Compute the first part of the local expansion}. For each box in level 2 and 3, we translate the multipole expansions from boxes in its interaction list into local expansions by
\begin{equation}\label{eq:compute local expansion (1)}
    \Psi^{b,(1)}_k = \sum_{b' \in \mathcal{L}_b} \sum_{k'=1}^{p} \mathbf{T}^{M2L}(k,k')\Phi^{b'}_{k'},
\end{equation}
where $b=4,\cdots,83$, $k=1,\cdots,p$ and $b'$ is in the interaction list of $b$.

\noindent \textbf{Step 4:} \textit{Compute the second part of the local expansion}. We transfer the local expansions for each box in level 2 to their children in level 3 using the $\mathbf{T}^{L2L}$ translation operator:
\begin{equation}\label{eq:compute local expansion (2)}
    \Psi^{b,(2)}_k = \sum_{k'=1}^{p} \mathbf{T}^{L2L}(k,k')\Psi^{b',(1)}_{k'},
\end{equation}
where $b=20,\cdots,83$, $k=1,\cdots,p$ and $b'$ is $b$'s parent.

\noindent \textbf{Step 5:} \textit{Complete the local expansion for the leaf boxes}. For each leaf box, the final local expansion is given by adding two parts together:
\begin{equation}\label{eq:complete local expansion}
    \Psi^{b}_k = \Psi^{b,(1)}_k + \Psi^{b,(2)}_k,
\end{equation}
where $b=20,\cdots,83$ and $k=1,\cdots,p$.

\noindent \textbf{Final Step:} \textit{Adding the far-field calculation to the near-field calculation}. The final result for each leaf box ($b = 20,\cdots, 83$) is obtained by adding the far-field calculation to the near-field calculation:
\begin{equation} \label{eq:final step}
    q_i = q_i^{(1)} + q_i^{(2)},
\end{equation}
where $i \in \mathbf{I}_b$. The near-field calculation is given by 
\begin{equation} \label{eq:final step: near-field calculation}
    q_i^{(1)} = \sum_{j \in \mathbf{I}_{\mathcal{N}_b}}a_{i,j}d_j.
\end{equation}
The far-field calculation is given by 
\begin{equation} \label{eq:final step: far-field calculation}
    q_i^{(2)}  \approx \sum_{k=1}^{p} u^{tgt,b}_{k,i}s^{tgt,b}_{k}\Psi^b_k.
\end{equation}

We note that the stepwise algorithm presented does not give a necessary order for the calculations. For example, the computation of local expansions for leaf boxes in Step 3 can be carried out once Step 1 is done.  

\subsection{Accuracy}\label{sec:accuracy}

In this section we discuss which factors affect the accuracy of the SVD-FMM. The first aspect to consider is the relative magnitude of the near-field and far-field terms in Eq. \eqref{eq:near and far-field mat-vec product box}. The near-field term is computed directly, while the far-field term is approximated. Numerical roundoff error for the direct calculation of the near field term is expected to be small. Therefore, the main source of error for the SVD-FMM is from the far-field calculation. Thus, if the magnitude of the near-field term is much larger than the far-field term, the error in the final result due to the SVD-FMM approximation will tend to be small. 

The accuracy of the far-field calculation is determined by the accuracy of the truncated SVD of the sub-matrices $\mathbf{A}(\mathbf{I}_{\mathcal{F}_b},\mathbf{I}_b)$ or $\mathbf{A}(\mathbf{I}_b, \mathbf{I}_{\mathcal{F}_b})$ in Eq. \eqref{eq:decomposition} and the accuracy of the translation operators given by Eqs. (\ref{eq:M2M} - \ref{eq:L2L}). The truncation error in Eq. \eqref{eq:decomposition} is dependent on the number of singular vectors used ($p$) and the $(p+1)$th singular value of the sub-matrices $\mathbf{A}(\mathbf{I}_{\mathcal{F}_b},\mathbf{I}_b)$ or $\mathbf{A}(\mathbf{I}_b, \mathbf{I}_{\mathcal{F}_b})$ \citep[e.g.,][Fact 9.14.28]{MatrixMathematics}. Additionally, \cite{gimbutas2003generalized} show that the error bounds on the translation operators rely on the $(p+1)$th singular value of the sub-matrices. Thus, the optimal number of the singular vectors used in the SVD-FMM may depend on the application. It should be determined by considering the trade-off between accuracy and efficiency, as the more the singular vectors are used, the more accurate the results but the slower the computation. 

The SVD-FMM algorithm uses truncated SVDs (Eq. \eqref{eq:decomposition}) and thus, may cause rank deficiency of the matrix $\mathbf{A}$. Nevertheless, the Hessian of the full variational problem will still be full rank \citep[e.g.,][]{tabeart2021new}. In addition, the matrix $\mathbf{A}$ used in the solution of the variational problem is not exactly the inverse of the matrix $\mathbf{R}$. However, only the matrix $\mathbf{A}$ (and not $\mathbf{R}$) is used in the solution of the variational problem. Furthermore, studies have shown that even approximate forms of spatial observation error correlations provide significant benefits to analysis accuracy compared with diagonal approximations \citep[e.g.,][]{StewartEtAl2008,stewart2013data,HealyAndWhite2005}.

\subsection{Algorithmic Complexity}\label{sec:complexity}

In this section we compare the number of floating point operations (flops) required for computing matrix-vector products using the standard approach and the SVD-FMM. Let $s=\overline{m_b}$ be the average number of observations in a leaf box (for $b$ in the highest level), and $B=m/s$ be the number of leaf boxes, where $m$ is the total number of observations. Computing Eq. \eqref{eq:compute multipole childless level} for a fixed $b$ and a fixed $k$ requires $2m_b$ operations ($m_b$ multiplication and $m_b$ add operations). Hence, calculating Eq. \eqref{eq:compute multipole childless level} for each value of $k$ requires $2m_bp$ operations. Finally, for each value of $b$, Eq. \eqref{eq:compute multipole childless level} requires $\sum_{b=20}^{83}2m_bp = 2mp$ operations. Similarly, Eq. \eqref{eq:compute multipole parent level} requires $2Bp^2$ operations. Eq. \eqref{eq:compute local expansion (1)} requires less than $2 \times 27Bp^2$ operations for leaf boxes and $2 \times 12 \times B/4 \times p^2$ operations for boxes in level 2, because there are at most 27 entries in the interaction list of each leaf box and 12 of each box at level 2. Eq. \eqref{eq:compute local expansion (2)} requires $2Bp^2$ operations. Eq. \eqref{eq:final step: near-field calculation} requires at most $2 \times 9ms$ operations, since the maximum number of boxes in the near field of each box is 9. Finally, Eq. \eqref{eq:final step: far-field calculation} requires $2mp$ operations. Thus, the operation count of the entire algorithm (excluding the initialization step) approximately sums to $18ms+4mp+64Bp^2$ or $18ms+4mp+64(m/s)p^2$.

The SVDs and translation operators only need to be computed once if the distribution of observations does not change. The singular vectors, singular values and translation operators can be used for any vector $\mathbf{d}$. However, if the distribution of the observations varies too often, then the computational costs in the initialization step and the computational cost of inverting the observation error covariance matrix could make the algorithm very expensive. 

For comparison, the direct matrix-vector multiplication of Eq. \eqref{eq:mat-vec product} requires $2m^2$ operations. In data assimilation, $\mathbf{q}$ is often computed using a forward and backward substitution \citep{golub1996cf,weston2014accounting,simonin2019pragmatic} which also requires $\mathcal{O}(m^2)$ operations. Fig. \ref{fig:complexity} compares the operation counts for direct computation and the SVD-FMM with our configuration of boxes and $p=10$. We observe that once the number of observations exceeds 500, the SVD-FMM requires fewer floating point operations than direct matrix-vector multiplication, and the difference increases with the number of observations.

\subsection{Parallelization}\label{sec:parallelism}

In this section we describe a novel parallel algorithm for the SVD-FMM. The FMM has several opportunities for parallelism \citep{greengard1990parallel} and our approach is not the only possibility. We have presented this section in order to provide a preliminary exploration of the potential of the SVD-FMM as a practical technique for operational data assimilation. However, we note that for our experimental results (section \ref{sec:results}) we have not used this parallel algorithm. The number of observations considered in our idealized experiments is much smaller than in operational applications, so that the serial calculations can be done within an acceptable time.

We should first distribute the matrix and vector elements across PEs according to the partitioning of observations in the box-tree. For each leaf box $b$, we should assign the sub-vector of $\mathbf{d}$ given by $\mathbf{d}(\mathbf{I}_{b})$ and the sub-matrix of $\mathbf{A}$ given by $\mathbf{A}(\mathbf{I}_b, \mathbf{I}_b \cup \mathbf{I}_{\mathcal{F}_b})$ to one PE. For each box $b$ in level 2, we should allocate the sub-matrix of $\mathbf{A}$ given by $\mathbf{A}(\mathbf{I}_b,\mathbf{I}_{\mathcal{F}_b})$ to one PE. For our particular configuration of boxes, we have 64 leaf boxes and 16 non-leaf boxes in level 2, hence we could use 80 PEs. However, to make our parallelization scheme easily comparable with the parallel formulations of matrix-vector multiplication given in sections \ref{sec:row-wise}-\ref{sec:symmetric}, we actually discuss the case where we choose 16 PEs out of 64 and let them store the data for level 2 boxes. We describe a possible parallelization of each mathematical step of the SVD-FMM as follows:

\begin{itemize}
\item \textit{Parallelization of the Initialization Step}. 
\begin{itemize}
  \item Each PE can calculate Eq. \eqref{eq:decomposition} independently. Once the singular vectors and singular values are obtained, the sub-matrices stored on each PE can be discarded.
  \item To compute Eq. \eqref{eq:M2M}, each PE that is assigned to a parent box should send $\mathbf{v}_{k}^{src}$ to three PEs that are assigned to its children. The calculation of the M2M translation operators would require 16 one-to-all broadcast operations to be performed simultaneously. The message size for each operation should be $mp/16$ and $4$ PEs should be involved.
  \item To compute Eq. \eqref{eq:M2L}, each PE should send a portion of $\mathbf{v}_{k}^{tgt}$ to at most 27 PEs, because the maximum number of boxes in an interaction list is 27. The calculation of the M2L translation operators should require an all-to-all broadcast with the message size of $mp/B$.
  \item To compute Eq. \eqref{eq:L2L}, each PE should send a portion of $\mathbf{v}_{k}^{tgt}$ to the PE that is assigned to its parent. The computation of the L2L translation operators would require 16 all-to-one reduction operations to be carried out in the same time. The message size for each operation should be $m_{\mathcal{F}_b}p$ for $b$ being a level 2 box. The communications required for calculating three translation operators can be done together using an all-to-all broadcast. After the initialization step, each PE that is assigned only to a leaf box will store $\mathbf{v}^{src}_k$, $s_k$, $\mathbf{T}^{M2M}$ and $\mathbf{T}^{M2L}$, and those assigned to both a leaf box and a level 2 box will store $\mathbf{T}^{M2L}$ and $\mathbf{T}^{L2L}$.
\end{itemize}
\item \textit{Parallelization of Step 1}. This step is perfectly parallel.
\item \textit{Parallelization of Step 2}. Each PE should compute $\mathbf{T}^{M2M}\cdot \Phi$ and then send the result to the PE that is assigned to its parent. Computing the multipole expansions for level 2 boxes would require an all-to-one reduction. The message size for this communication operation should be $p$.
\item \textit{Parallelization of Step 3}. Each PE should compute $\mathbf{T}^{M2L}\cdot \Phi$ and then send the result to another PE. Each PE should collect the partial result of Eq. \eqref{eq:compute local expansion (1)} from at most 27 PEs. This step would need an all-to-all broadcast with a message size of $p$ or $2p$.
\item \textit{Parallelization of Step 4}. Each PE that is assigned to a box in level 2 should compute $\mathbf{T}^{L2L}\cdot \Phi$ and then send the result to the PEs that are assigned to its children. This step would use an one-to-all broadcast with a message size of $p$.
\item \textit{Parallelization of Step 5}. This step is perfectly parallel. 
\item \textit{Parallelization of final step}. The far-field calculation for each leaf-box is perfectly parallel. For the computation of the near-field calculation, each PE should obtain the elements of $\mathbf{d}$ from at most 8 PEs, which is the maximum number of one box's neighbours. This would require an all-to-all broadcast. The average message size for this operation is $m/B$. We can use one all-to-all broadcast to complete the communication tasks for this step and step 3.
\end{itemize}

Table \ref{table} summarises the communication costs for the SVD-FMM and four parallel formulations of the matrix-vector multiplications described in section \ref{sec:parallel standard mat-vec multiplication}. The major advantage of using the SVD-FMM is that the message size for each communication operation is dramatically reduced.

In the proposed parallelization scheme, we suggested using $B$ PEs and letting one of every four PEs conduct the computations for boxes in level 2. This could be particularly useful if the supercomputer configuration has non-uniform memory access (NUMA) nodes with locally shared memory \citep{grama2003introduction}. In this case, we could avoid the remote communications required in Step 2, Step 4 and similar steps in the Initialization. Alternatively, we could assign the tasks for level 2 and level 3 to different PEs, i.e., using 80 PEs for our configuration of boxes. This could allow each PE to carry out a similar amount of work. More generally, we note that the practical implementation of the given parallelization scheme should be adjusted to suit the configuration of the supercomputer. 

\begin{center}
\captionof{table}{The communication time for four distinct parallel formulations of matrix-vector multiplication with full matrices and the parallelization scheme for the SVD-FMM. The number of the observations is $m$, the number of leaf boxes is $B$ and the number of singular vectors is $p$. The \# PEs represents the number of PEs that are involved in a communication operation. The symbols $t_s$ and $t_w$ are the startup time and the per-word transfer time respectively, which are determined by the configuration of the parallel machine.}
\begin{tabular}{c|c|c|c|c}\label{table}
 \textbf{Parallelization Scheme} & \textbf{Communication Operation} & \textbf{\# PEs} & \textbf{Message Size}& \textbf{Communication Time} \\ 
 \hline
 Row-wise & all-to-all broadcast & $B$ & $m/B$ &$t_s\log B+t_wm$\\ 
 \hline
 Column-wise & all-to-one reduction &$B$ & $m$&$(t_s+t_wm)\log B$ \\
 & scatter&$B$&$m/B$&$t_s\log B+t_wm$ \\ 
 \hline
 Block & one-to-all broadcast &$\sqrt{B}-1$&$m/\sqrt{B}$&$(t_s+t_wm/\sqrt{B})\log\sqrt{B}$ \\ 
 & all-to-one reduction &$\sqrt{B}-1$ &$m/\sqrt{B}$&$(t_s+t_wm/\sqrt{B})\log\sqrt{B}$ \\ 
 \hline
 Symmetric & all-to-all broadcast & $B$ & $\approx m/B$ & $< t_s\log B+t_wm$\\ 
 & all-to-all reduction & $B$ & $\approx m/B$ & $< t_s\log B+t_wm$ \\
 \hline
 SVD-FMM & all-to-one reduction &4&$p$& $(t_s+t_wp)\log(4)$\\ 
 & all-to-all broadcast &$B$& $p$, $2p$ or $m/B$ & $< t_s\log B+t_wm$\\ 
 & one-to-all broadcast &4 & $p$ & $(t_s+t_wp)\log (4)$\\
 \hline
\end{tabular}
\end{center}

\section{Experimental Design}\label{sec:experimental design}

Our numerical experiments are designed to demonstrate the potential of applying the SVD-FMM to compute the matrix-vector products involved in operational data assimilation. We investigate the accuracy of using the SVD-FMM under different scenarios that may occur in practical applications. We conduct serial calculations rather than using parallel computing because in this initial study we have chosen to study idealized problems of moderate size as this is sufficient for our goals.

For this initial study, we have not included matrix inversion as part of the SVD-FMM algorithm. However, to obtain our experimental results we require an inverse observation error covariance matrix. We use the INV function in MATLAB \citep{inv} to compute the inverses of $\mathbf{R}$ and reconditioned $\mathbf{R}$. The generation of the $\mathbf{R}$ matrix is described in section \ref{sec:Modelling of the observation error covariance matrix} and reconditioning methods are given in section \ref{sec:Reconditioning of the observation error covariance matrix}. The INV function uses an LDL decomposition \citep{golub1996cf}. The required SVDs of sub-matrices of $\mathbf{R}^{-1}$ (denoted $\mathbf{A}$) are computed using the SVDS function in MATLAB \citep{svds}, which uses a Lanczos method and is efficient for finding a few singular values and vectors of a large matrix \citep{Larsen_1998,baglama2005augmented}. 

In the results section (section \ref{sec:results}) we use the log of root-mean-squared error (RMSE) to assess the accuracy of SVD-FMM, which is defined as
\begin{equation}\label{eq:log(RMSE)}
   \log(\text{RMSE})=\log( \sqrt{\frac{\sum_{i=1}^{m}(q^{fmm}_i-q^{ref}_i)^2}{m}}),
\end{equation}
where superscript $fmm$ denotes the matrix-vector product computed using the SVD-FMM and $ref$ denotes the matrix-vector product obtained by the standard approach. Note that the $\log(\text{RMSE})$ shown in section \ref{sec:results} is averaged over 100 realizations of $q^{fmm}_i$, where each one uses a different $\mathbf{d}$.

\subsection{Observation distribution and observation-minus-model departures}

To simulate observation data, we assume our observations are regularly distributed over a region from $54^{\circ}$ to $60^{\circ}$N and $6^{\circ}$W to $6^{\circ}$E with a grid-length of 12 km. This is similar to moderately thinned geostationary satellite data used in operational forecasting over the UK \citep{waller2016SEVIRI}. This results in 3456 observations and an error covariance matrix of size $3456 \times 3456$. For convenience, we directly sample the observation-minus-model departure vector $\mathbf{d}$ from a Gaussian distribution with mean zero and (innovation) covariance given by $\mathbf{R}$ plus background error covariance. The background error covariance is modelled using the SOAR correlation function with a lengthscale of 20 km and a standard deviation of 0.6 (section \ref{sec:experimental design}). The correlation lengthscale and standard deviation are selected according to \cite{Ballard2016} to be appropriate for km-scale numerical weather prediction. The results presented in section \ref{sec:results} have been averaged over 100 realizations of $\mathbf{d}$.

In practice, the observation distribution may be different at each assimilation cycle because of factors such as quality control. Therefore, in some of our experiments we have also applied the SVD-FMM to compute the matrix-vector product with randomly chosen missing observations.

\subsection{Modelling of the observation error covariance matrix} \label{sec:Modelling of the observation error covariance matrix}

Modelling observation error covariance matrices using correlation functions is an useful approach to deal with observation error correlations \citep[e.g.,][]{stewart2013data,tabeart2018conditioning}. Here we use several correlation functions to model the observation error covariance matrices. The first is the Gaussian correlation function \citep[e.g.,][]{haben2011conditioning},
\begin{equation}
    \mathbf{C}(i,j)=\exp(-\frac{|r_{i,j}|^2}{2l^2}),
\end{equation}
where $l>0$ is the correlation length scale and $r_{i,j}$ denotes the great-circle distance between two observations. The second is the first-order auto-regressive (FOAR) correlation function, also called the Markov correlation function \citep[e.g.,][]{stewart2013data},
\begin{equation}
    \mathbf{C}(i,j)=\exp(\frac{-|r_{i,j}|}{l}).
\end{equation}
We also use the SOAR correlation function \citep[e.g.,][]{daley1994atmospheric,tabeart2018conditioning},
\begin{equation}
    \mathbf{C}(i,j)=(1+\frac{|r_{i,j}|
    }{l})\exp(\frac{-|r_{i,j}|}{l}),
\end{equation}
and the Mat\'ern 5/2 correlation function, \citep[e.g.,][]{rasmussen2006gaussian}, 
\begin{equation}
    \mathbf{C}(i,j)=(1+\frac{\sqrt{5}r_{i,j}}{l}+\frac{5r_{i,j}^2}{3l^2})\exp(-\frac{\sqrt{5}r_{i,j}}{l}).
\end{equation}
The observation error covariance matrix can be generated using the correlation functions by
\begin{equation}
    \mathbf{R}=\mathbf{D}\mathbf{C}\mathbf{D},
\end{equation}
where $\mathbf{D}$ is a diagonal matrix whose diagonal elements are a prescribed standard deviation. For our experiments, we choose the standard deviation as one and the correlation lengthscales as $l=80$, $160$, $240$ km. These correlation lengthscales are selected based on the horizontal error correlations estimated for geostationary satellite data used in operational km-scale data assimilation \citep{waller2016SEVIRI} and for AMVs \citep{Cordoba17}. It should be noted that the Markov and Mat\'ern 5/2 correlation functions may lead to a sparse $\mathbf{A}$, in which case the SVD-FMM may not be optimal to use. We use these correlation functions in our experiments because we wish to show that the SVD-FMM can be applied to a variety of matrices.

\subsection{Reconditioning of the observation error covariance matrix} \label{sec:Reconditioning of the observation error covariance matrix}

In practical applications, observation error covariance matrices are often ill-conditioned \citep{tabeart2020improving,Haben2011a}. Moreover, if the matrix $\mathbf{R}$ has a large condition number, this can lead to poor convergence of the minimisation problem in variational data assimilation \citep{weston2014accounting, tabeart2018conditioning,tabeart2020impact}. In our applications, the condition number of the matrices that are created using correlation functions is dependent on the chosen function and correlation lengthscale \citep{haben2011conditioning,Haben2011b}.

We use two common matrix reconditioning techniques to reduce the matrix condition number in our experiments: the ridge regression method and the minimum eigenvalue method \citep{tabeart2020improving,weston2014accounting,campbell2017accounting,bormann2016enhancing}. In the ridge regression (RR) method the reconditioned matrix is given by
\begin{equation}
    \mathbf{R}_{RR}=\mathbf{R}+\delta\mathbf{I},
\end{equation}
where $\mathbf{I}$ is the identity matrix and 
\begin{equation}\label{eq:RR}
    \delta=\frac{\lambda_{max}-\lambda_{min}\kappa_{req}}{\kappa_{req}-1},
\end{equation}
where $\lambda_{max}$ is the maximum eigenvalue of $\mathbf{R}$, $\lambda_{min}$ is the minimum eigenvalue of \textbf{R} and $\kappa_{req}$ is the required new condition number. 

The minimum eigenvalue (ME) method changes the eigenvalue spectrum of the matrix $\mathbf{R}$ by setting all eigenvalues smaller than a threshold value to the threshold value. The threshold value ($T$) is given in terms of the required condition number as
\begin{equation}\label{eq:ME-threshold}
    T=\lambda_{max}/\kappa_{req}.
\end{equation}
The reconditioned matrix is then constructed via 
\begin{equation}
    \mathbf{R}_{ME}=\mathbf{E}\Lambda_{ME}\mathbf{E}^\mathrm{T},
\end{equation}
where $\mathbf{E}$ is a square matrix whose columns are the eigenvectors of $\mathbf{R}$ and $\Lambda_{ME}$ is a diagonal matrix with the diagonal elements being the new eigenvalues. 

\section{Results}\label{sec:results}

\subsection{The number of singular vectors and the size of singular values}

\begin{figure}[!ht]
	\centering
	\includegraphics{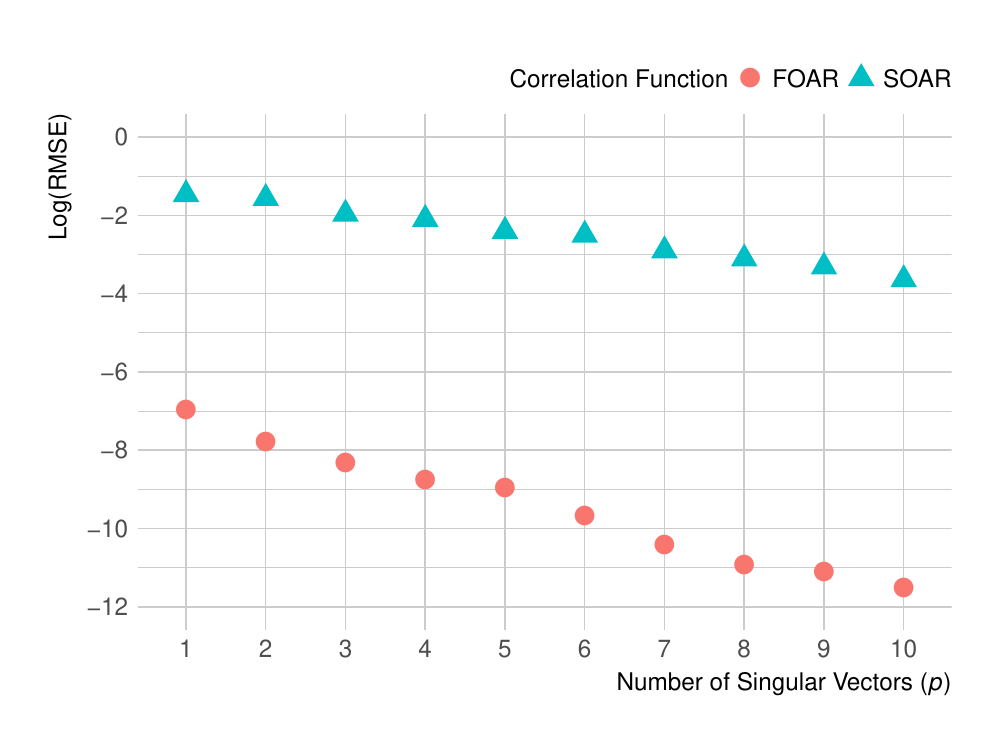}
	\caption{The $\log(\text{RMSE})$ for the SVD-FMM against the number of singular vectors $(p)$ used in the approximation. The matrices $\mathbf{A}$ are given by inverting the FOAR and SOAR covariance matrices with correlation lengthscale $l=80$ km. The $\log(\text{RMSE})$ is averaged over 100 realizations of $\mathbf{d}$.}
	\label{fig:sing_vec}
\end{figure}

\begin{figure}[!ht]
	\centering
	\includegraphics{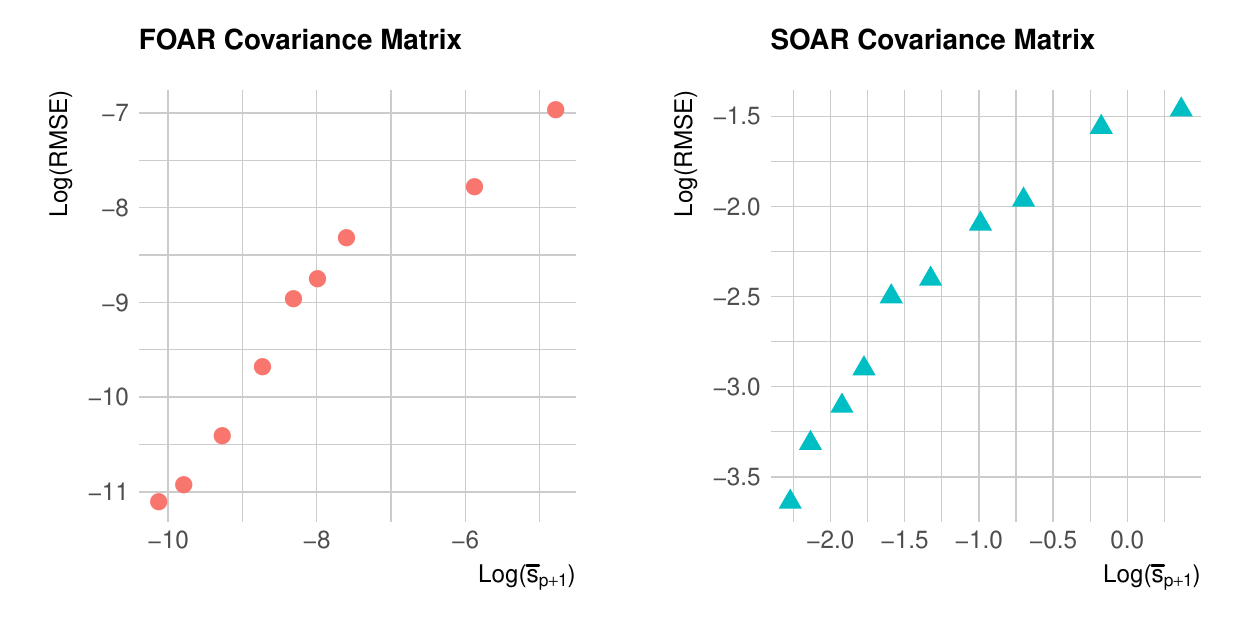}
	\caption{The $\log(\text{RMSE})$ for the SVD-FMM using $p$ singular vectors against the log of mean ($p+1$)th singular value of the sub-matrices of the inverse FOAR and SOAR covariance matrices used in Fig. 8.}
	\label{fig:sing_val}
\end{figure}

We carried out an experiment to assess the accuracy of the SVD-FMM, as the number of singular vectors ($p$) changes. As discussed in section \ref{sec:accuracy}, we expect the accuracy to depend on both $p$ and the $(p+1)$th singular values for each sub-matrix. Fig. \ref{fig:sing_vec} provides the $\log(\text{RMSE})$ for the SVD-FMM with the FOAR and SOAR covariance matrices as a function of $p$. As anticipated, the $\log(\text{RMSE})$ decreases, and hence the accuracy increases, as more singular vectors are used. Moreover, Fig. \ref{fig:sing_val} demonstrates that the $\log(\text{RMSE})$ for the SVD-FMM using $p$ singular vectors has an approximately linear relationship with the log of the mean ($p+1$)th singular value of the sub-matrices, which is given by 
\begin{equation}
    \log(\overline{s}_{p+1})=\log(\overline{s}_{p+1}^{src,b})=\log(\overline{s}_{p+1}^{tgt,b})
\end{equation}
for $p=1, \cdots, 9$ and $b=4, \cdots, 83$. Additionally, we note that the SVD-FMM with the FOAR covariance matrix is more accurate than with the SOAR covariance matrix for all the values of $p$ considered. This is because the mean $(p+1)$th singular value of the sub-matrices of the FOAR matrix is consistently smaller than that of the SOAR matrix (Fig. \ref{fig:sing_val}). We have also carried out several other experiments using different correlation lengthscales, different correlation functions and reconditioned covariance matrices and found similar qualitative results: the accuracy of the SVD-FMM using $p$ singular vectors relies on the $(p+1)$th singular values of the sub-matrices. 

These results using the mean singular values of the sub-matrices provide some guidance as to how to set the value of $p$ in a given application. However, as they require the computation of the SVDs of all of the relevant sub-matrices of $\mathbf{A}$, they are time-consuming to compute. 
\subsection{Varying correlation lengthscale}

\begin{figure}[!ht]
	\centering
	\includegraphics{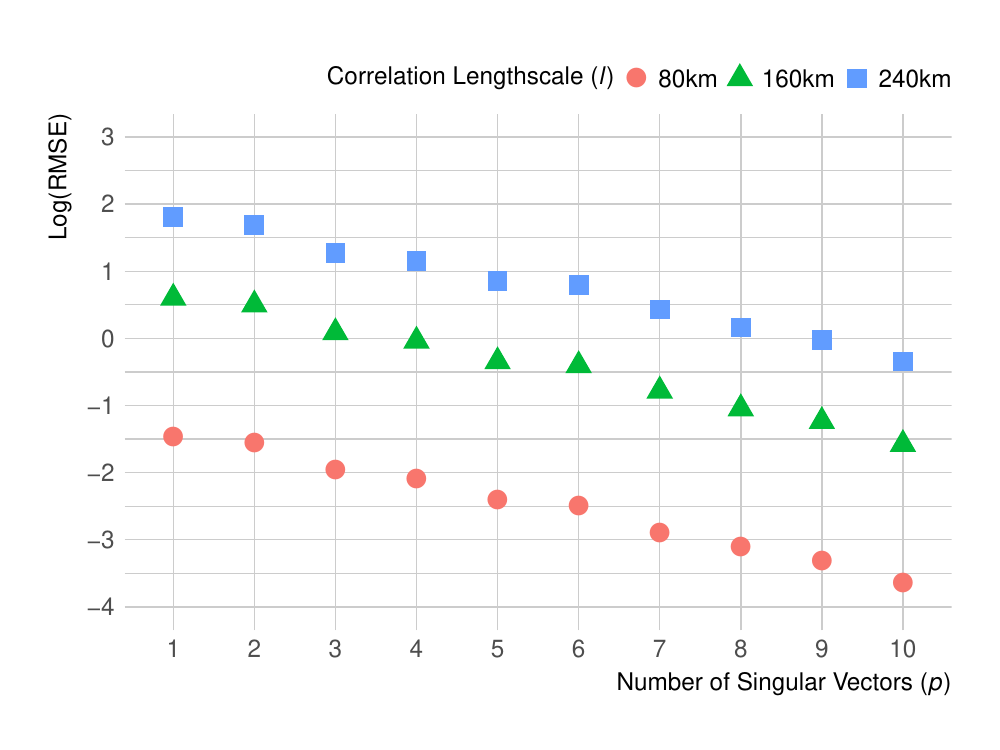}
	\caption{As Fig. 8, but matrices $\mathbf{A}$ are given by inverting the SOAR covariance matrices with three correlation lengthscales.}
	\label{fig:lengthscale}
\end{figure}

Different observations can exhibit correlated errors over different lengthscales. To examine how correlation length affects the accuracy of the SVD-FMM, we use three correlation lengthscales for the SOAR correlation function: $l=80$, $160$ and $240$ km. Fig. \ref{fig:lengthscale} reveals an increase of the $\log(\text{RMSE})$ of the SVD-FMM with correlation lengthscale. Nevertheless, we still only need a few singular vectors to obtain a small $\log(\text{RMSE})$. However, if we desire to obtain a given accuracy, we may need to use a larger value of $p$ for longer lengthscales. In addition, we note that the magnitude of the elements of $\mathbf{A}$ increases as correlation lengthscale increases. For a smaller lengthscale, the far-field elements of $\mathbf{A}$ are close to zero and could be neglected to give a sparse matrix. Moreover, compared to the SOAR correlation function, the Mat\'ern and Markov functions are more likely to give a sparse $\mathbf{A}$ if the correlation lengthscale is small.

The variation of the accuracy of the SVD-FMM with correlation lengthscale can also be explained by the singular values of the sub-matrices. We find in our numerical experiments (not shown) that the singular values of the sub-matrices are smaller when the correlation lengthscale is smaller. It is known that the maximum singular value of the full matrix (inverse SOAR covariance matrix) also decreases as the correlation lengthscale reduces \cite[section 5.3.2]{haben2011conditioning}. Therefore, we investigated whether it would be possible to use the singular values of the full matrix to estimate the accuracy of the SVD-FMM. \citet{THOMPSON19721} showed that the singular values of the sub-matrices are bounded by the singular values of the full matrix, i.e.,
\begin{equation}\label{eq:inequality}
    \sigma_i(\mathbf{B}) \le \sigma_i(\mathbf{A}), \; i=1, 2, 3, \ldots, \min \{ \alpha, \beta \}
\end{equation}
where $\mathbf{B} \in \reals^{ \alpha \times \beta}$ denotes a sub-matrix of $\mathbf{A}$, $\sigma_i(\mathbf{B})$ denotes the $i$-th singular value of $\mathbf{B}$ and $\sigma_i(\mathbf{A})$ denotes the $i$-th singular value of $\mathbf{A}$. Equality can be obtained for some choices of $\mathbf{B}$. However, since the dimensions of the sub-matrices used in SVD-FMM are much smaller than the dimensions of the full matrix, our numerical experiments showed that $\sigma_i(\mathbf{B})$ is typically much smaller than $\sigma_i(\mathbf{A})$ in practice. We compared the maximum singular values of different full matrices numerically and found that they could provide a rough guide to relative accuracies: for a given value of $p$, the best accuracy was obtained for the full matrix with the smallest maximum singular value. 

\subsection{Reconditioned covariance matrices}

\begin{figure}[!ht]
	\centering
	\includegraphics{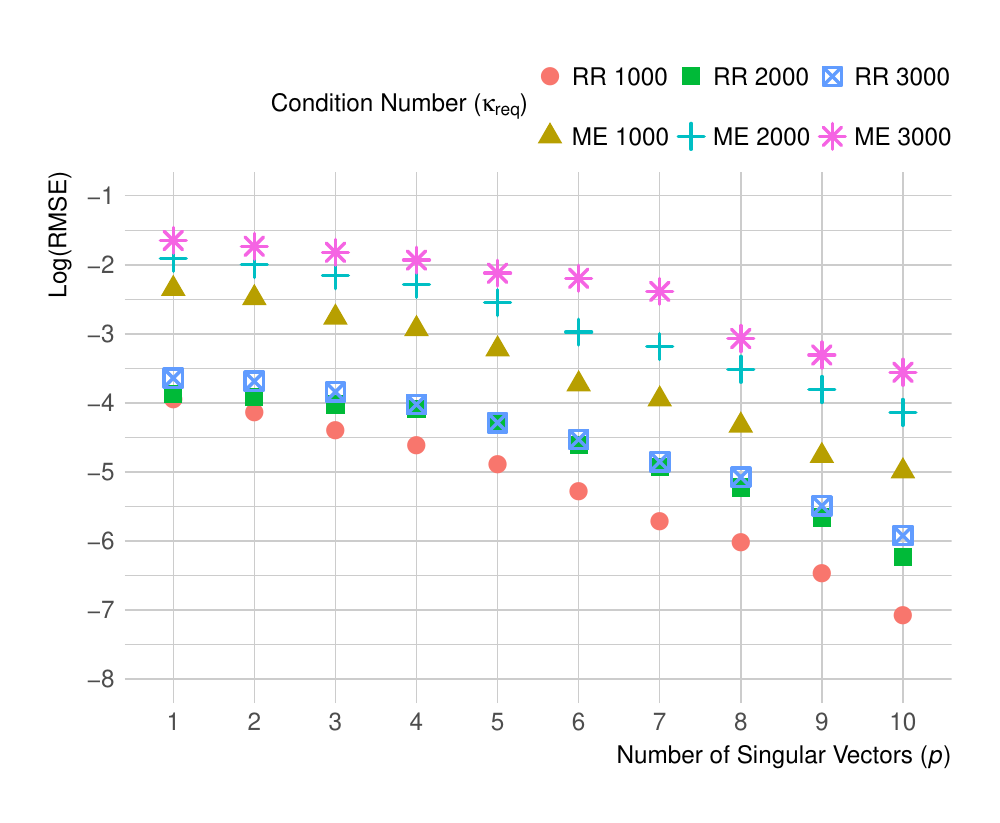}
	\caption{As Fig. 8, but matrices $\mathbf{A}$ are given by inverting reconditioned SOAR covariance matrices with correlation lengthscale $l=80$ km. The ridge regression (RR) and minimum eigenvalue (ME) methods are used to reduce the condition number of SOAR covariance matrix to target values ($\kappa_{req}$).}
	\label{fig:reconditioning}
\end{figure}

The observation error covariance matrices used in data assimilation procedures are often ill-conditioned and thus, in order to invert them we need to reduce their condition numbers using proper techniques \citep{tabeart2020improving}. Reconditioning the covariance matrices may also improve the convergence behaviour of the minimisation procedure \citep{weston2014accounting,tabeart2020impact,tabeart2021new}. We evaluate how reconditioning affects the accuracy of the SVD-FMM; we reduce the condition number of the SOAR covariance matrix with $l=80$ km to three required new condition numbers ($\kappa_{req}=1000$, $2000$ and $3000$) using the RR and ME methods. Fig. \ref{fig:reconditioning} shows that reconditioning the SOAR correlation matrix can improve the accuracy of SVD-FMM; the smaller the condition number of the SOAR covariance matrices after reconditioning, the better the accuracy. 

We also find that the RR method gives smaller $\log(\text{RMSE})$ than the ME method for the same required condition number. This difference is caused by the different ways the two reconditioning methods reduce the condition number, which leads to different size singular values of the full inverse covariance matrices satisfying:
\begin{equation}\label{eq:singular value RR ME}
    \sigma_{max}(\mathbf{A}_{RR}) < \sigma_{max}(\mathbf{A}_{ME}),
\end{equation}
where $\sigma_{max}(\mathbf{A}_{RR})$ and $\sigma_{max}(\mathbf{A}_{ME})$ denote the maximum singular value of the reconditioned matrices obtained using the RR and ME methods, respectively. Combing Eq. \eqref{eq:inequality} and Eq. \eqref{eq:singular value RR ME}, the leading singular values of the sub-matrices of $\mathbf{A}_{RR}$ are expected to typically be smaller than the leading singular values of the sub-matrices of $\mathbf{A}_{ME}$. Hence, the SVD-FMM with reconditioned matrices using RR method will usually have smaller error. A derivation of Eq. \eqref{eq:singular value RR ME} is presented in Appendix.

\subsection{Varying correlation functions}

\begin{figure}[!ht]
	\centering
	\includegraphics{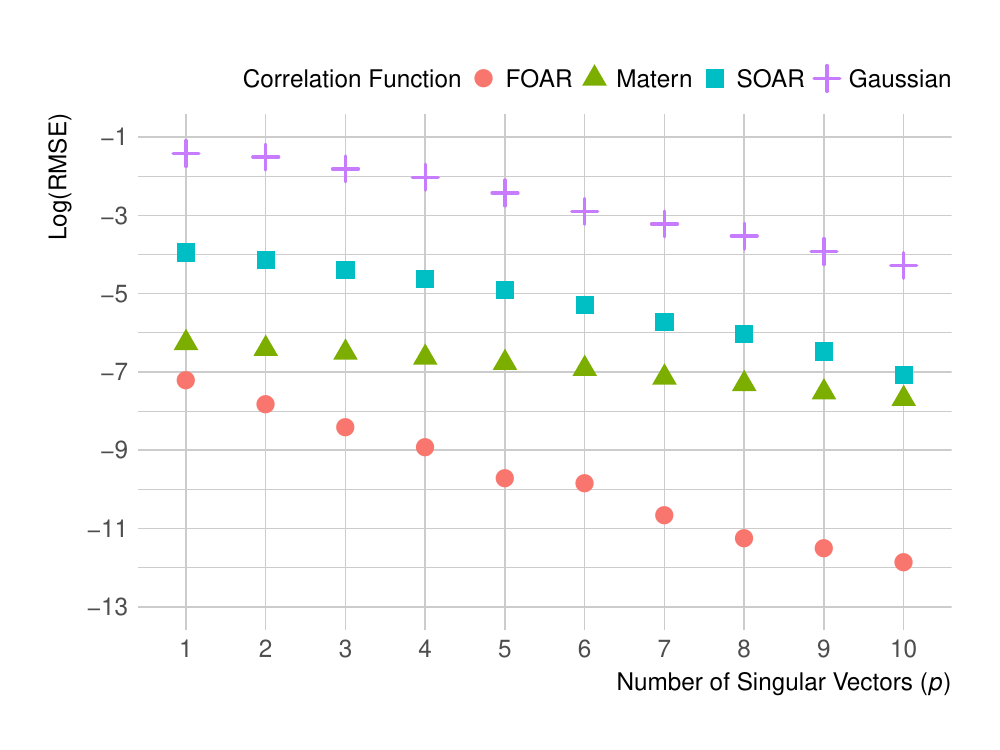}
	\caption{As Fig. 8, but matrices $\mathbf{A}$ are given by inverting reconditioned Gaussian, FOAR, SOAR and Matern covariance matrices with a correlation lengthscale of 80 km. All covariance matrices are reconditioned to have a new condition number of 1000 using the RR method prior to inversion.}
	\label{fig:correlation function}
\end{figure}

To provide more numerical evidence for the wide applicability of the SVD-FMM, we use other correlation functions than SOAR and FOAR to create covariance matrices. To allow for comparison with Fig. \ref{fig:reconditioning} we use the RR method to reduce the condition number of all the matrices to 1000 before inverting them. Fig. \ref{fig:correlation function} shows that the SVD-FMM can work well with the inverses of the matrices given by different correlation functions. We note that although the condition number of each matrix is identical after reconditioning, the accuracy of the SVD-FMM still relates to the original condition numbers; the larger the condition number prior to reconditioning, the greater the $\log(\text{RMSE})$. In our experiments the FOAR correlation matrix has a condition number on the order of thousands, while the Gaussian correlation matrix is extremely badly conditioned and has a condition number on the order of $10^{15}$ prior to reconditioning. 

By comparing Fig. \ref{fig:correlation function} with Fig. \ref{fig:reconditioning}, we also observe that reducing the condition number of the Gaussian correlation matrix to 1000 using the RR method gives larger $\log(\text{RMSE})$ than reducing the condition number of the SOAR correlation matrix to 3000 using the same method. Therefore, to acquire a comparable accuracy using the SVD-FMM with the same value of $p$, we may need to reduce the condition number of two matrices to different values. The matrix with a larger condition number prior to reconditioning may need a smaller condition number after reconditioning.  

\subsection{Removing a portion of observations}

\begin{figure}[!ht]
	\centering  
	\includegraphics{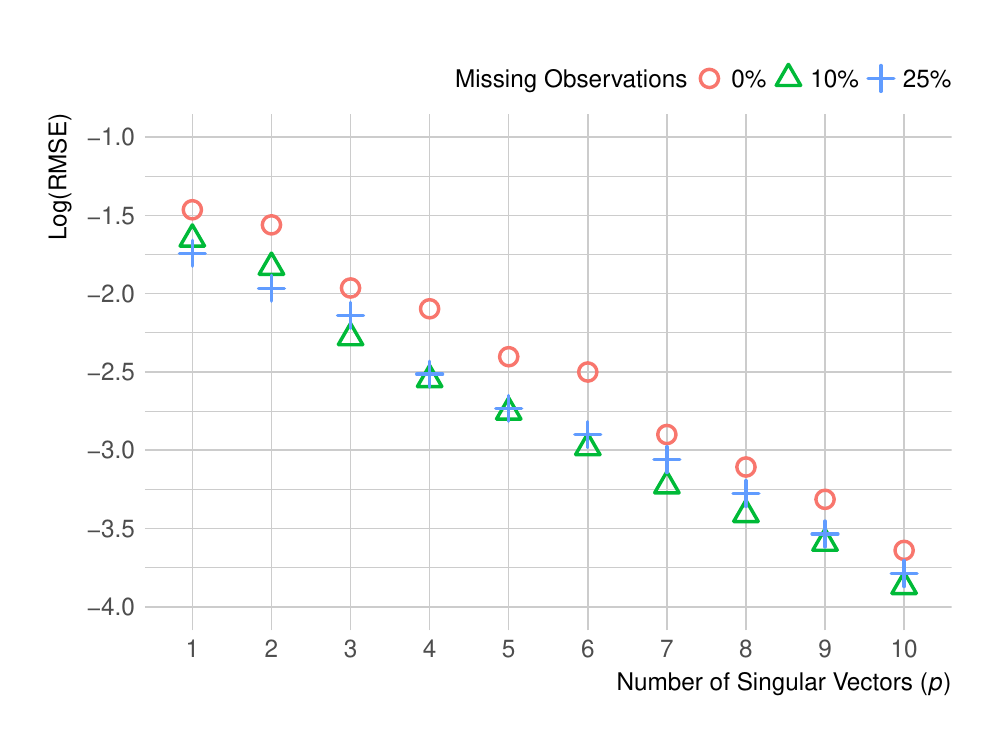}
	\caption{As Fig. 8, but matrices $\mathbf{A}$ are given by inverse SOAR covariance matrices with correlation lengthscale $l=80$ km and randomly chosen missing observations.}
	\label{fig:missing observation}
\end{figure}

In operational data assimilation, the number of observations varies each assimilation cycle due to quality control, among other reasons. This will lead to a decrease of the dimension of covariance matrix and disrupt the structure of the matrix. The resultant covariance matrix lacks some of the rows and columns of the original one. The missing observations may cause a problem if a leaf box becomes empty or contains fewer observations than $p$. This could be solved by using a different partition of observations. In our experiment, we randomly select the locations of the missing observations, and with up to 25\% missing observations, every leaf box always contains enough observations. Fig. \ref{fig:missing observation} shows that the SVD-FMM can still work well with missing observations without losing accuracy. In our experiments using an inverse SOAR covariance matrix, we find that the missing observations actually lead to a slight decrease in the $\log(\text{RMSE})$ compared to the full set of observations. The difference in the log(RMSE) between removing 10\% of the observations and removing 25\% of the observations is not statistically significant.

\section{Conclusion and discussion}\label{sec:summary}

Some observations have been shown to exhibit strong spatial error correlations, e.g., Doppler radar radial winds, geostationary satellite data, and AMVs. Accounting for this information in data assimilation systems can improve analysis accuracy and forecast skill. However, it can make the computation of the products of observation weighting matrices and observation-minus-model departure vectors very expensive, in terms of not only the computational complexity, but also the communication costs in parallel computing. We have therefore investigated the use of the SVD-FMM for the rapid computation of these matrix-vector products. This numerical approximation method is best suited for full or dense precision matrices. If the observation weighting matrix is sparse it might be appropriate to only compute the near-field calculation. Moreover, the SVD-FMM is suitable for use for large problems, when the overhead of computing the SVD-FMM expansions is outweighed by the reductions in floating-point operations and communication costs compared with alternative approaches, and where each PE has enough memory to store its part of the error covariance matrix of these observations. The exact range of the number of observations that makes the SVD-FMM useful depends on the computer architecture and the constraints of the operational schedule.

We proposed a novel possible parallelization scheme for the SVD-FMM in our application. In comparison to the parallel formulations of direct matrix-vector multiplication (section \ref{sec:parallel standard mat-vec multiplication}), the parallelization scheme for the SVD-FMM largely reduces the size of the messages transferred. This reduces communication costs, making the use of full observation error covariance matrices associated with large spatial extents potentially feasible in operational data assimilation. 

In our idealised experiments we have examined the accuracy of the SVD-FMM, in terms of applying it to the inverses of various observation error covariance matrices. These matrices are created using commonly used correlation functions, such as Gaussian and SOAR correlation functions, and different correlation lengthscales. In some of our experiments, we have applied reconditioning methods to the covariance matrices before inverting them. We have also investigated the feasibility of the SVD-FMM with randomly chosen missing observations. We find consistent results as discussed in section \ref{sec:accuracy}: i) the accuracy of the SVD-FMM increases as more singular vectors are used and ii) the variation of the accuracy is approximately linear with the mean spectrum of singular values of the sub-matrices used in the approximation. Furthermore, our experiments indicate that a comparison of the maximum singular values of the full inverse matrices themselves can be used as a rough guide to determine which matrices will give more accurate results with the SVD-FMM.

The most computationally expensive parts of our current implementation of the SVD-FMM are the inversions of covariance matrices and the SVDs of the sub-matrices of inverse matrices. They need to be re-performed every time the observation error covariance matrix is changed. This happens when the number of observations or the distribution of observations are changed. Nevertheless, our work has shown that the SVD-FMM has potential for use in operational data assimilation for fast computation of the products of the inverse observation error covariance matrices and observation-minus-model departure vectors. This will allow a large volume of observational data to be assimilated within a short time interval, which is particularly important for applications such as convection-permitting hazardous weather forecasting.

\section*{ACKNOWLEDGEMENTS}
The authors would like to express their thanks to Amos Lawless, Nancy K. Nichols, David Simonin and Joanne A. Waller for their useful discussions. G. Hu and S. L. Dance were funded in part by the UK EPSRC DARE project (EP/P002331/1). The MATLAB code for SVD-FMM developed by the authors is available on request from the corresponding author.  

\section*{CONFLICT OF INTEREST}
We declare we have no competing interests.

\bibliographystyle{abbrvnat}
\bibliography{ref}

\section*{APPENDIX: MAXIMUM SINGULAR VALUES OF RECONDITIONED MATRICES USING RR AND ME METHODS}\label{appendix:RR ME}

In this appendix we derive the result given in Eq. \eqref{eq:singular value RR ME}. The RR method increases each of eigenvalues of the original matrix by the same amount via adding an increment given by Eq. \eqref{eq:RR} to the diagonal. In contrast, the ME method changes only the smallest eigenvalues to a value given by Eq. \eqref{eq:ME-threshold}. Let $\lambda_{min}(\mathbf{R}_{RR})$ denote the minimum eigenvalue of $\mathbf{R}_{RR}$ and $\lambda_{min}(\mathbf{R}_{ME})$ denote the minimum eigenvalue of $\mathbf{R}_{ME}$, then we have
\begin{equation}
    \lambda_{min}(\mathbf{R}_{RR})=\frac{\lambda_{max}-\lambda_{min}\kappa_{req}}{\kappa_{req}-1}+\lambda_{min} > \frac{\lambda_{max}-\lambda_{min}\kappa_{req}}{\kappa_{req}}+\lambda_{min}=\frac{\lambda_{max}}{\kappa_{req}} = \lambda_{min}(\mathbf{R}_{ME}).
\end{equation}
The reciprocal of the minimum eigenvalue of a matrix is the maximum eigenvalue of its inverse and hence we have 
\begin{equation}
    \lambda_{max}(\mathbf{A}_{RR})=\frac{1}{\lambda_{min}(\mathbf{R}_{RR})} < \frac{1}{\lambda_{min}(\mathbf{R}_{ME})}=\lambda_{max}(\mathbf{A}_{ME}).
\end{equation}
Since the eigenvalues and singular values of a symmetric positive semidefinite matrix are the same \citep[][Definition 5.6.1]{MatrixMathematics}, we obtain
\begin{equation}
    \sigma_{max}(\mathbf{A}_{RR}) = \lambda_{max}(\mathbf{A}_{RR}) < \lambda_{max}(\mathbf{A}_{ME}) = \sigma_{max}(\mathbf{A}_{ME})
\end{equation}
as required.

\end{document}